\documentclass[preprint,12pt]{elsarticle}

\usepackage{amssymb}
\usepackage{amsmath}
\usepackage{graphicx}
\usepackage{caption}
\usepackage{subcaption}
\usepackage{amssymb}
\usepackage{amsmath}
\usepackage{array}
\usepackage{mathtools}

\usepackage{lineno}

\journal{Computers and Geotechnics}

\begin{document}
\graphicspath{FIGS/}
\begin{frontmatter}



\title{Development of a stable two-phase contact MPM algorithm for saturated soil-structure interaction problems}


\author [label1]{Chihun Sung}
\author [label2]{Shyamini Kularathna}
\author [label1]{Krishna Kumar}

\affiliation[lable1]{organization={Department of Civil, Architectural, and Environmental Engineering, University of Texas at Austin},
            city={Austin}, 
            state={Texas},
            postcode={78712},
            country={United States}}
            
\affiliation[label2]{organization={Department of civil and environmental engineering, University of California, Berkeley, California, 94720, United States}
}
             
\begin{abstract}
The simulation of soil-structure interaction problems involving two-phase materials poses significant challenges in geotechnical engineering. These challenges arise due to differences in material stiffnesses, interaction between multiple phases, high bulk modulus of pore fluid, and low permeability. The conventional explicit time integration scheme is limited by its conditional stability, necessitating small time step sizes and resulting in pressure oscillations under rapid loading conditions. To address these issues, we propose a stable two-phase contact algorithm within the Material Point Method (MPM) framework for soil-structure interaction problems. Our algorithm models the soil as a fully saturated porous media with incompressible pore fluid. We introduce three main advancements over conventional MPM methods. We employ Chorin's projection method to solve coupled formulations and reduce numerical oscillations. By implicitly handling a diffusion term, our algorithm permits larger stable time step sizes, independent of the bulk modulus and permeability of the pore fluid. Lastly, We integrate a rigid algorithm to model solid bodies accurately and a precise contact detection algorithm. We provide detailed formulations and time increment processes of the two-phase contact MPM algorithm. Furthermore, we compare the proposed algorithm with Finite Element Method (FEM) and explicit MPM to assess its accuracy and performance in simulating coupled hydro-mechanical problems. The two-phase contact algorithm offers a more stable and efficient approach to simulate soil-structure interaction problems.
\end{abstract}


\end{frontmatter}


\section{Introduction}
\label{S:1}
Soil-structure interaction is the reaction between soil and adjacent structures under varying loading conditions. One key task of the soil-structure interaction problem is estimating the deformation of the structures on top of or embedded in soil under various loading conditions. These loading conditions may include direct loadings on the structures or wave loadings provoked by an earthquake or moving traffic. Predicting saturated soil reaction to a complex and dynamic loading from the upper structure is another integral part of the soil-structure interaction problems. With recent developments in powerful computers and numerical tools, soil-structure interaction problems can now be solved numerically, leading to rapid advancements in the field. \citep{Kausel2010}.

However, there remain challenges in simulating geotechnical problems involving soil-structure interaction. In particular, when modeling soil as a two-phase material body in problems such as pipeline embedment, footing penetration, and slope stability problems, a contact model must be able to handle the interaction between material bodies with significantly different stiffnesses and different phases. Furthermore, numerical simulations are conditionally stable due to the high bulk modulus and low permeability of pore fluid. In addition, as these problems often involve large displacements resulting in significant geometry changes, a proper numerical tool or technique is necessary. Given these challenges, numerical modeling of saturated soil-structure interaction problems with large displacement has been a research subject of great interest in recent geotechnical engineering.

The finite element method (FEM), a classical mesh-based numerical method, has been extensively utilized for early applications of multiphase formulations, as shown in \citep{ghaboussi1972variational}, \citep{zienkiewicz1977coupling}, \citep{zienkiewicz1984dynamic}, \citep{prevost1984localization}, \citep{simon1992multiphase}, \citep{borja1998elastoplastic}. As an accurate and robust tool, FEM is suitable for modeling continuum bodies in a wide range of mechanical applications. However, when dealing with simulations involving large deformations, severe mesh distortion issues inevitably arise. To address these challenges, complex remeshing techniques such as Arbitrary Lagrangian-Eulerian (ALE) and Coupled Eulerian-Lagrangian (CEL) were introduced to the FEM, as discussed in \citep{hirt1974arbitrary}, \citep{hu1998practical}, \citep{benson1992computational}. \citep{hamann2015application} applied the CEL method to model pile penetration into the saturated soil. \citep{wang2015large} introduced three ALE-related methods to handle large-displacement geotechnical problems. \citep{souli2001arbitrary} used the ALE method to model free surface fluid flow problems. Despite their effectiveness in resolving mesh distortion issues, these methods can be computationally costly and yield inaccurate state variables when modeling history-dependent materials.

The material point method (MPM) has emerged as a powerful tool in the geotechnical field for simulating various soil mechanics problems with large deformations. Developed by \citep{Sulsky1995}, MPM is a continuum particle-based method that uses an Eulerian-Lagrangian spatial discretization. It discretizes a deformable material body into Lagrangian material points and solves the governing equations on the Eulerian background grid. Thanks to its spatial discretization scheme and continuum framework, we can efficiently model the macroscopic multi-phase granular problems with large deformations. Moreover, since the Lagrangian material points carry the history-dependent state variables, it can efficiently model history-dependent materials. Applications of the MPM for hydro-mechanical interface problems include \citep{Zambrano-Cruzatty2020}, \citep{Zhang2009}, \citep{Phuong2016}. For more detailed review of large-deformation MPM simulations, \citep{Soga2016} can be referred.

However, the conventional two-phase MPM applications have limitations as they typically use the explicit time integration scheme. Numerical simulations involving explicit time integrations for two-phase problems are only conditionally stable, and the critical time step size is limited by the compressibility and permeability of the pore fluid. This limitation results in an unreasonably small time step size. Additionally, when the pore fluid is nearly incompressible or incompressible, explicit simulations suffer from volume-locking problems, which yield inaccurate patterned pore pressure distributions. To ensure stability and circumvent these problems, a mixed-order spatial integration scheme for pressure and displacement must be used (\citep{Gresho1998}). However, even after meeting these stability conditions, the inherent instability of the explicit time integration scheme can still produce non-physical oscillations in the pressure field (\citep{Bandara2015}).

Another main drawback of the conventional MPM applications for soil-structure problems is the inaccurate contact detection process. The multi-velocity field contact algorithm proposed by \citep{Bardenhagen2000} identifies the interface nodes when nearby particles of different materials are approaching each other. This algorithm identifies the contact significantly earlier than the material bodies are actually in contact, thus generating an artificial gap between the bodies. The lack of a softening function also causes stress oscillation; the contact algorithm with an explicit time integration scheme often produces artificial stress or pressure variations at the contact surface (\citep{ma2014new}). Finally, because the MPM reconstructs the background grid to its original position at the end of the time step, the impenetrability condition is not satisfied at the beginning of each time step. This condition would contribute to inaccurate initial velocities of the material bodies and eventually inaccurate results (\citep{huang2011contact}). Although using denser mesh can reduce the errors, it also increases computational costs exponentially. Therefore, this study aims to propose a two-phase stable contact MPM algorithm for saturated soil-structure interaction problems.

To achieve this objective, the study adopts a semi-implicit MPM algorithm that was previously developed for a saturated porous media by \citep{kularathna2018splitting} for a saturated porous media and extend this algorithm by introducing contact and a rigid body algorithm. The core attributes of the semi-implicit MPM algorithm are: (1) an introduction of Chorin’s projection method (\citep{chorin1968numerical}), and (2) an implicit approach to the diffusion equation using an iterative solver. The splitting algorithm provides numerical stability and computational efficiency by dividing a time step into two sub-steps and applying equal-order approximation functions to pressure and displacement variables. The semi-implicit Euler time integration scheme can minimize the limitations of the explicit scheme as it significantly reduces the artificial pressure oscillation and allows larger time steps independent of the fluid incompressibility and permeability constraints (\citep{Kularathna2017}).

The treatment of the contact algorithm in this study ensures accuracy and stability, thus considerably reducing the problems of the conventional multi-velocity field MPM contact algorithm. By introducing the projection method, the proposed algorithm adjusts the material velocities at the interface three times in a time step to satisfy the impenetrability condition: at the beginning of a time step, between the sub-steps, and at the end of the time step.olid bodies are modeled as rigid, allowing the use of a time increment step independent of material stiffness and eliminating contact problems between bodies with significantly different stiffnesses. Additionally, we introduced the Generalized Interpolation Material Point (GIMP) method (\citep{bardenhagen2004generalized}) to avoid the cell-crossing error of the linear MPM and to take particle size into account in the interface detection process. In this contact algorithm, we track the distances between particle edges near a background node, and the node is correctly identified to be at the interface when the material bodies are in contact. 
This paper present the detailed semi-implicit contact MPM formulation and a rigid body algorithm. The performance of the proposed algorithm is evaluated by comparing its results with the analytical solution and those from the conventional FEM and the explicit MPM, using example simulations of interface problems.

\section{Material Point Method}
\label{S:2}
The MPM is a particle-based Eulerian-Lagrangian computational method where a continuum body is discretized into Lagrangian material points. The material points can move through the background Eulerian grid, carrying the physical and history-dependent properties of the continuum. The governing equations are solved at the background grid nodes, but the position of the nodes is not updated during the computation. After solving the governing equation, the kinematic variables are mapped from the nodes to material points with the shape functions.

\begin{figure}[h]
    \centering
    \includegraphics[width=0.8\textwidth]{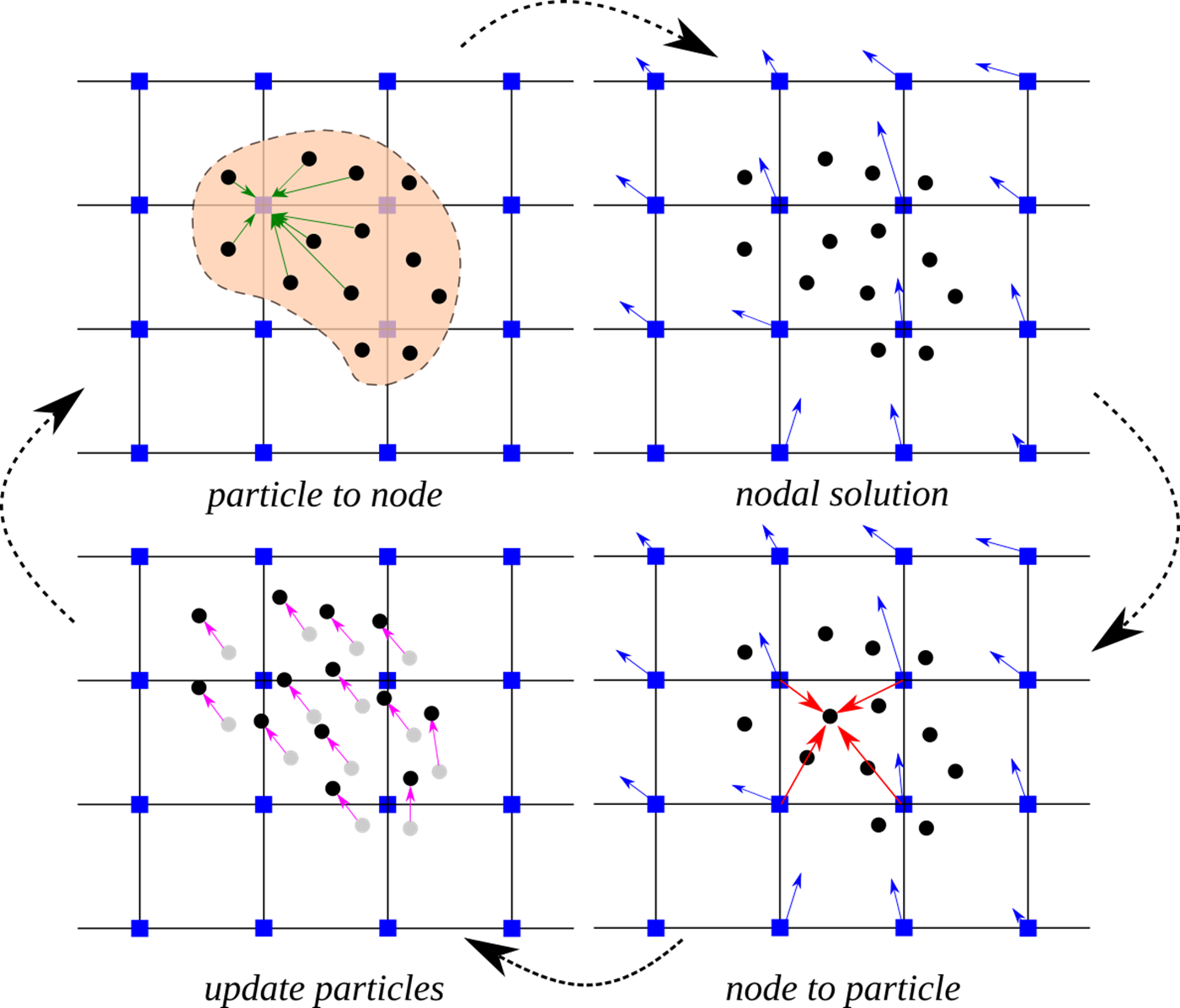}
    \caption{Illustration of the single-point two phase MPM model}
    \label{fig:Illustration of the MPM algorithm}
\end{figure}

Figure \ref {fig:Illustration of the MPM algorithm} illustrates the MPM algorithm. In the initial phase (particle to node), the properties carried by material points are projected to the adjacent nodes (velocity, volume, mass, and forces). In phase 2 (nodal solution), the governing equations of motion are solved at the nodes, and nodal velocities and positions are computed. In phase 3 (node to particle), the computed nodal kinematics are mapped back to the material points. Finally, in the last phase (update particles), the positions of the material points are updated. 

In the MPM, non-slip contact is naturally achieved between material bodies as their kinematics belong to a single velocity field. To simulate the relative motion at the interface, \citep{Bardenhagen2000} proposed a multi-velocity approach in which different nodal velocity fields are utilized for different material bodies. The nodes at the interface are detected by comparing the velocity of each body to the combined center-of-mass velocity.

The multi-velocity contact algorithm is briefly presented below. The pre-designation of the interface is not required as the algorithm identifies the interface during the simulation. Nodes at the interface of more than one material body are identified by comparing the velocity of each body to the center-of-mass velocity:

\begin{equation}
    \mathbf{v}^a \neq \mathbf{v}^{CoM}
    \label{eq:contact_1}
\end{equation}

Where the superscript  is the index of the material bodies and $\mathbf{v}^{CoM}$ is the center-of-mass velocity. The surface unit normal vector is used to identify the approaching and separation of the bodies. A body is approaching another material body when:

\begin{equation}
    (\mathbf{v}^{CoM} - \mathbf{v}^a) \cdot \mathbf{n}^a > 0
    \label{eq:contact_2}
\end{equation}

Where $\mathbf{n}$ is the unit surface normal vector for each material body. When the approach is identified, the contact formulation is enforced to the interface node if the normal surface traction is compressive, and free separation is allowed otherwise. To prevent interpenetration, the velocities of each material body in the normal direction are adjusted to the center-of-mass velocity: 

\begin{equation}
    \Delta \mathbf{v}^\alpha_{normal}=(\mathbf{v}^{CoM} - \mathbf{v}^a) \cdot \mathbf{n}^a
    \label{eq:contact_3}
\end{equation}

\begin{equation}
    \mathbf{v}^a=\mathbf{v}^a+\Delta \mathbf{v}^ \alpha _{normal} \cdot \mathbf{n}^ \alpha
    \label{eq:contact_4}
\end{equation}
Frictional contact can be applied in the tangential direction. This approach allows separation, sliding, and rolling while forbidding inter-penetration.

\section{Semi-implicit two-phase MPM formulation}
\label{S:3}

\citep{kularathna2018splitting} implemented the semi-implicit two-phase algorithm in the framework of the single-point two-phase MPM, where each material point consists of a solid skeleton and pore fluid. The porous medium is assumed to be completely saturated with incompressible water. We adopt the $\mathbf{u}-p$ formulation thus ignoring the Darcy velocity. The u-p formulation is convenient for static and quasi-static simulations but suffers from an accuracy issue in modeling problems with high-frequency loadings. We present the detailed formulation of the semi-implicit two-phase model.

The conservation of mass for the mixture considering incompressible water and solid is:

\begin{equation}
    n\nabla\cdot\mathbf{v}_w+(1-n) + (1-n)\nabla\cdot\mathbf{v}_s=(\nabla\cdot\mathbf{v}_s+n\nabla\cdot\mathbf{v}_{wR})=0
    \label{eq:two-phase_1}
\end{equation}
Where $n$ is the porosity, $\mathbf{v}_s$ is the solid velocity, $mathbf{v}_w$ is the water velocity, and $\mathbf{v}_{wR}=\mathbf{v}_w - \mathbf{v}_s$ is the seepage velocity. The dynamic momentum conservation for the mixture is:

\begin{equation}
    \rho^s \frac{\partial\mathbf{v}_s}{\partial t}+\rho^w \frac{\partial\mathbf{v}_w}{\partial t} = \nabla\cdot(\mathbf{T}^s_E - p\mathbf{I}) +\rho g
    \label{eq:two-phase_2}
\end{equation}
Where $rho$ is the mass density of the mixture, $\mathbf{T}^s_E$ is the effective stress tensor, $p$ is the pore pressure, and $\mathbf{I}$ is the identity tensor. $\rho^w$ and $\rho^s$ are the mass density of water and solid, respectively:
\begin{equation}
    \rho^w = n \rho^{wR},\hspace{12pt} \rho^s = (1-n) \rho^{sR}
    \label{eq:two-phase_3}
\end{equation}
Where $rho^{sR}$ is the mass density of solid grain, $rho^{sR}$ is the mass density of water grain. The momentum balance of the water is:
\begin{equation}
    \rho ^w \frac{\partial\mathbf{v}_w}{\partial t}=\nabla\cdot (-np\mathbf{I})+\rho^w\mathbf{g}+\hat{\mathbf{p}}^F
    \label{eq:two-phase_4}
\end{equation}
Where $\hat{\mathbf{p}}^F$ is the drag force, defined as:
\begin{equation}
    \hat{\mathbf{p}}^F=-\frac{n^2\rho ^{wR} \mathbf{g}}{k}\mathbf{v}_{wR}+p\nabla n
    \label{eq:two-phase_5}
\end{equation}
Where $k$ is the permeability of pore fluid. The time-discretized formulations conservation of mass for the mixture phase is:
\begin{equation}
    n\nabla\cdot\mathbf{v}^{t+1}_w+(1-n)\nabla\cdot\mathbf{v}^{t+1}_s=0
    \label{eq:time-two-phase_1}
\end{equation}
The time-discretized momentum balance equations for the mixture phase and water phase are:
\begin{equation}
    \rho^s\frac{\mathbf{v}^{t+1}_s-\mathbf{v}^{t}_s}{\Delta t}+\rho^w\frac{\mathbf{v}^{t+1}_w-\mathbf{v}^{t}_w}{\Delta t}=\nabla\cdot(\mathbf{T}^s_E)-\nabla p^{t+1}+\rho\mathbf{g}
    \label{eq:time-two-phase_2}
\end{equation}
\begin{equation}
    \rho^w\frac{\mathbf{v}^{t+1}_w-\mathbf{v}^{t}_w}{\Delta t}=-n\nabla p^{t+1}+\rho^w\mathbf{g}-\frac{n^2p^{wR}g}{k}\mathbf{v}_{wR}
    \label{eq:time-two-phase_3}
\end{equation}
Where $\Delta t$ is the time step size. 

The time-discretized formulations are divided into two sub-steps using Chorin's projection method. In the first sub-step, the intermediate velocities of the solid and water are computed, while the pore pressure is handled explicitly. By introducing the intermediate velocities, the conservation of mass for the mixture is expressed as:
\begin{equation}
    n\nabla\cdot\mathbf{v}^*_w+(1-n)\nabla\cdot\mathbf{v}^*_s=\Delta t (\frac{1-n}{\rho^{sR}}+\frac{n}{\rho^{wR}})\nabla^2(p^{t+1}-p^t)
    \label{eq:chorin_conservation_mass}
\end{equation}
Where superscript $*$ denotes the intermediate parameters. $\mathbf{v}^*_s$ is the intermediate solid velocity and $\mathbf{v}^*_w$ is the intermediate water velocity. The momentum balance equations for the mixture phase using intermediate velocities are computed as:
\begin{equation}
    \rho^s\frac{\mathbf{v}^{*}_s-\mathbf{v}^{t}_s}{\Delta t}+\rho^w\frac{\mathbf{v}^{*}_w-\mathbf{v}^{t}_w}{\Delta t}=\nabla\cdot(\mathbf{T}^s_E)-\nabla p^{t}+\rho\mathbf{g}
    \label{eq:chorin_momentum_balance_mixture1}
\end{equation}
\begin{equation}
    \rho^s\frac{\mathbf{v}^{t+1}_s-\mathbf{v}^{*}_s}{\Delta t}+\rho^w\frac{\mathbf{v}^{t+1}_w-\mathbf{v}^{*}_w}{\Delta t}=-\nabla (p^{t+1}-p^{t+1})
    \label{eq:chorin_momentum_balance_mixture2}
\end{equation}
The momentum balance equations for the water phase are:
\begin{equation}
    \rho^w\frac{\mathbf{v}^{*}_w-\mathbf{v}^{t}_w}{\Delta t}=-n\nabla p^{t}+\rho^w\mathbf{g}-\frac{n^2\rho^{wR}g}{k}\mathbf{v}^*_wR
    \label{eq:chorin_momentum_balance_water1}
\end{equation}
\begin{equation}
    \rho^w\frac{\mathbf{v}^{t+1}_w-\mathbf{v}^{*}_w}{\Delta t}=-n\nabla (p^{t+1}-p^{t})-\frac{n^2\rho^{wR}g}{k}(\mathbf{v}_{wR}-\mathbf{v}^*_{wR})
    \label{eq:chorin_momentum_balance_water2}
\end{equation}
Where $\mathbf{v}^*_{wR}$ is the intermediate seepage velocity.

Applying the $\mathbf{u}-p$ formulation which ignores the seepage acceleration, the intermediate and final water accelerations are replaced with solid accelerations:
\begin{equation}
    \mathbf{a}^*_w = \frac{\mathbf{v}^*_w+\mathbf{v}^t_w}{\Delta t}=\mathbf{a}^*_s = \frac{\mathbf{v}^*_s+\mathbf{v}^t_s}{\Delta t}
    \label{eq:accleration_expressions1}
\end{equation}
\begin{equation}
    \mathbf{a}_w = \frac{\mathbf{v}^{t+1}_w+\mathbf{v}^t_w}{\Delta t}=\mathbf{a}_s = \frac{\mathbf{v}^{t+1}_s+\mathbf{v}^t_s}{\Delta t}
    \label{eq:accleration_expressions2}
\end{equation}
\begin{equation}
    \mathbf{a}_w - \mathbf{a}^*_w = \frac{\mathbf{v}^{t+1}_w+\mathbf{v}^*_w}{\Delta t}=\mathbf{a}_s-\mathbf{a}^*_s = \frac{\mathbf{v}^{t+1}_s+\mathbf{v}^*_s}{\Delta t}
    \label{eq:accleration_expressions3}
\end{equation}
Where $\mathbf{a}^*_w$ is the intermediate water acceleration, $\mathbf{a}_w$ is the final water acceleration, $\mathbf{a}^*_s$ is the intermediate solid acceleration, and $\mathbf{a}_w$ is the final solid acceleration. Finally, the mixture momentum balance equations are written as:
\begin{equation}
    \rho\mathbf{a}^*_s=\nabla\cdot(\mathbf{T}^E_s-p^t\mathbf{I})+\rho \mathbf{g}
    \label{eq:acc_chorin_mix_momentum1}
\end{equation}
\begin{equation}
    \rho(\mathbf{a}_s-\mathbf{a}^*_s)=-\nabla(p^{t+1}-p^t)
    \label{eq:acc_chorin_mix_momentum2}
\end{equation}
The water momentum balance equations are:
\begin{equation}
    \rho^w\mathbf{a}^*_s+\frac{n^2\gamma^{wR}}{k}\mathbf{v}^*_{wR}=-n\nabla p^t + \rho^w\mathbf{g}
    \label{eq:acc_chorin_water_momentum1}
\end{equation}
\begin{equation}
    \rho^w(\mathbf{a}_s-\mathbf{a}^*_s)+\frac{n^2\gamma^{wR}}{k}(\mathbf{v}_{wR}-\mathbf{v}^*_{wR})=-n\nabla (p^{t+1}-p^t)
    \label{eq:acc_chorin_water_momentum2}
\end{equation}
The mixture mass balance equation is:
\begin{equation}
    [\frac{k}{\rho^wg}(\frac{\rho^w}{\rho}-n)-\frac{\Delta t}{\rho}]\nabla^2 (p^{t+1}-p^t)=-\nabla\cdot\mathbf{v}^*_s-n\nabla\cdot\mathbf{v}^*_{wR}
    \label{eq:acc_chorin_mix_mass}
\end{equation}

The intermediate solid acceleration and intermediate water velocity are computed with the mixture and fluid momentum balance equations (equation \ref {eq:acc_chorin_mix_momentum1} and equation \ref {eq:acc_chorin_water_momentum1}), respectively. In the first sub-step, the intermediate velocities do not satisfy the incompressibility condition. After the pressure increment is updated with the mixture mass balance equation to satisfy the divergence of the flow velocity (equation \ref {eq:acc_chorin_mix_mass}), and finally, the intermediate solid acceleration is corrected to the final acceleration with the gradient of pressure increment (equation \ref{eq:acc_chorin_mix_momentum2}). 

\section{Semi-implicit Contact MPM Algorithm}
\label{S:4}

\begin{figure}[h]
    \centering
    \includegraphics[width=1\textwidth]{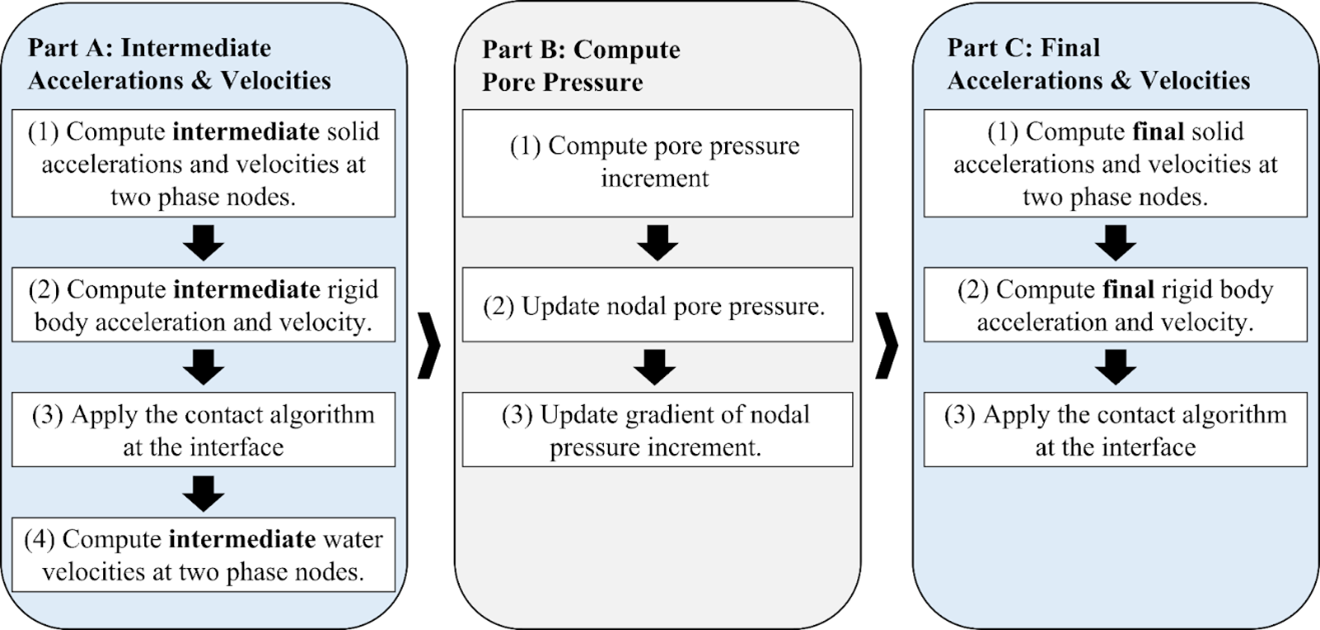}
    \caption{Flowchart of semi-implicit contact MPM algorithm solver}
    \label{fig:flowchart}
\end{figure}

In this chapter, we introduce a numerical algorithm for solving the semi-implicit two-phase MPM formulations in combination with the contact and rigid body algorithms. To illustrate the step-by-step procedure of each stage in the semi-implicit contact algorithm, we provide a flowchart in Figure \ref{fig:flowchart}. The algorithm comprises three primary stages: (a) computing intermediate acceleration and velocities, (b) updating pore pressure, and (c) computing final acceleration and velocities. We specifically consider the interface between a two-phase material body and a rigid body, where the two-phase material represents either porous media or soil. Finally, we present a contact detection algorithm that considers the particle sizes and locations.

\subsection{Numerically Discretized Formulations}
\subsubsection{Intermediate velocity and acceleration (step A)}
(a-1) Compute intermediate velocities of the two-phase material body

The nodal intermediate solid and water velocities are computed using the momentum balance equations. Nodal intermediate solid accelerations and velocities of the two-phase body are:
\begin{equation}
    \mathbf{a}^{s.tb^*}_i=\frac{\mathbf{b}^{mix,tb}_i+\mathbf{t}^{mix,tb}_i+\mathbf{f}^{mix,tb}_i}{m^{mix,tb}_i}
    \label{eq:int_tp_solid_acc}
\end{equation}

\begin{equation}
    \mathbf{v}^{s,tb^*}_i=\mathbf{v}^{s,tb}_{i,initial}+\Delta t \cdot \mathbf{a}^{s,tb^*}_i
    \label{eq:int_tp_solid_vel}
\end{equation}
Where, superscript $*$ denotes the intermediate kinematic parameters at the first sub time step, superscript $tb$ denotes the two-phase material body ($=tb$), superscript $mix$ indicates the mixture phase, and superscript $s$ indicates the solid phase. The subscript $i$ presents nodal parameters, and subscript $initial$ means initial kinematic parameters at the beginning of the time step. $\mathbf{a}^{s^*}_i$ is the nodal intermediate solid acceleration, $\mathbf{v}^{s^*}_i$ is the nodal intermediate solid velocity, and $\mathbf{v}_{i,initial}$ is the nodal initial solid velocity. $\mathbf{b}^{mix}_i$ is the nodal mixture body force, $\mathbf{t}^{mix}_i$ is the nodal mixture traction force, $\mathbf{f}^{mix,tb}_i$ is the nodal mixture internal force, $m^{mix}_i$ is the nodal mixture mass and $\Delta t$ is the time step size.

\hfill

\noindent(a-2) Compute intermediate rigid body acceleration\\

\begin{figure}[h]
    \centering
    \includegraphics[width=0.6\textwidth]{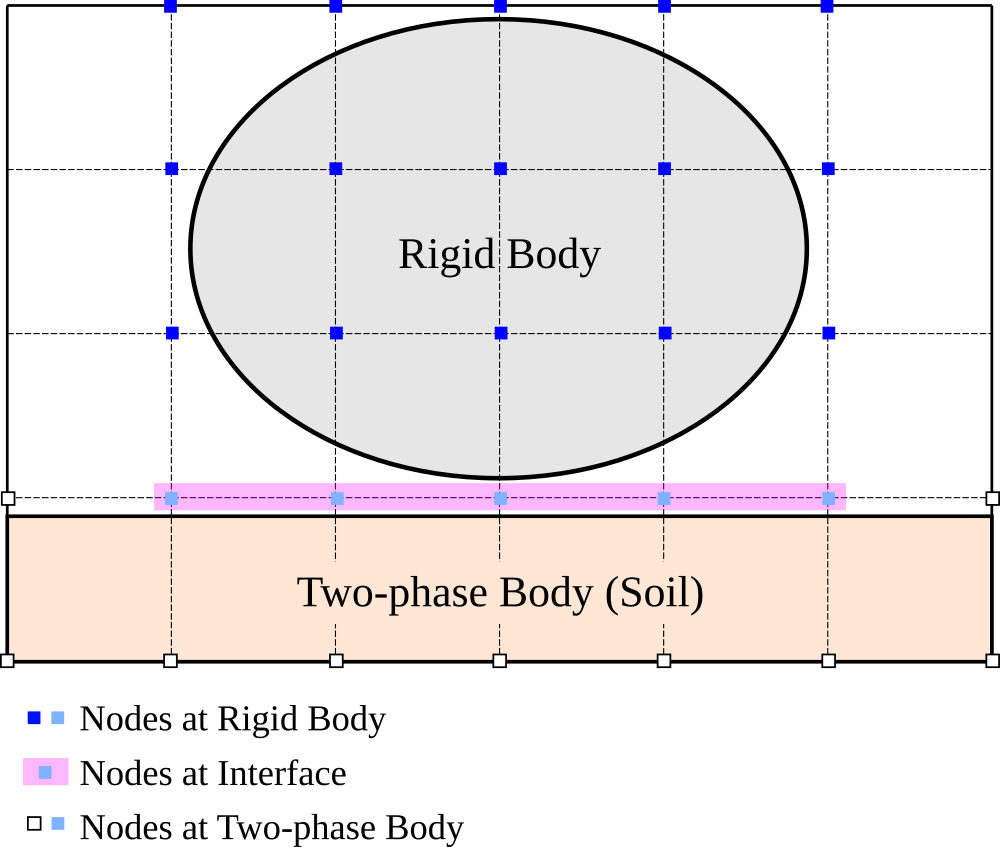}
    \caption{Illustration of the rigid body algorithm}
    \label{fig:rigid_algorithm}
\end{figure}

The rigid body acceleration is computed globally by combining nodal parameters belonging to the rigid body and those at the interface. Figure \ref{fig:rigid_algorithm} provides a visual representation of the nodes classification. Rigid body nodes refers to nodes adjacent to the rigid body particles, while the nodes at interface are the nodes adjacent to both the rigid body particles and two-phase material body particles. The intermediate rigid body acceleration, $\mathbf{a}^{rb^*}$, is calculated using a mixture momentum balance equation:
\begin{equation}
    \mathbf{a}^{rb^*}=\frac{\sum_{i=1}^{n_rb} (\mathbf{b}^{mix,rb}_i+\mathbf{t}^{mix,rb}_i)+\sum_{i=1}^{n_{contact}} (\mathbf{b}^{mix,tb}_i+\mathbf{t}^{mix,tb}_i)+\sum_{i=1}^{n_{contact}} \mathbf{f}^{mix,tb}_i}{\sum_{i=1}^{n_{contact}} m^{mix,CoM}_i}
    \label{eq:int_rigid_acc}
\end{equation}

Where $n_{rb}$ is the number of nodes belong to the rigid body, and $n_{contact}$ is the number of nodes that are identified as in contact, superscript $CoM$ represents the center-of-mass parameter, superscript $rb$ represents the rigid body parameter, and $m^{mix,CoM}_i$ is the nodal center-of-mass mixture mass. The numerator in equation \ref{eq:int_rigid_acc} has three terms. The first term represents the total mixture external force acting on the rigid body, which is the sum of the nodal mixture body forces and nodal traction forces that belong to the rigid body nodes. The second term is the total mixture external forces of the two-phase body at the interface. Finally, the third term is the total mixture internal force of the two-phase body at the interface. The computed intermediate rigid body acceleration is assigned to the nodes in the rigid body, which becomes the nodal intermediate solid velocity:
\begin{equation}
    \mathbf{a}^{s,rb^*}_i=\mathbf{a}^{rb^*}
\end{equation}

\hfill

\noindent(a-3) Apply contact algorithm at the interface\\

Applying the contact algorithm to the rigid body, the nodal intermediate center-of-mass solid velocity is set to the intermediate rigid body acceleration at the interface:
\begin{equation}
    \mathbf{a}^{s,CoM^*}_i=\mathbf{a}^{s,rb^*}_i
\end{equation}
We use the contact algorithm to adjust the nodal acceleration of the two-phase material body at the nodes at the interface. Unlike the conventional multi-velocity field contact algorithm, which adjusts only velocities, in the semi-implicit two-phase contact algorithm, the contact algorithm is applied to the solid acceleration in the first sub-step, as the intermediate water velocity is computed with the intermediate solid acceleration. As the rigid body in this study is frictionless, we only adjust the nodal component of the nodal intermediate solid acceleration. The subscript $normal$ indicates the surface normal component.
\begin{equation}
    \Delta\mathbf{a}^{s,tb^*}_{i,normal}=(\mathbf{a}^{s,CoM^*}_i-\mathbf{a}^{s,tb^*}_i)\cdot \mathbf{n}^{tb}_i
\end{equation}
\begin{equation} 
\mathbf{a}^{s,tb*}_i=\mathbf{a}^{s,tb*}_i+\Delta\mathbf{a}^{s,tb^*}_{i,normal}\cdot \mathbf{n}^{tb}_i
\end{equation}
The, the intermediate center-of-mass solid velocity at the interface is updated:
\begin{equation} 
\mathbf{v}^{s,tb*}_i=\mathbf{v}^{s}_{i,initial}+\mathbf{a}^{s,tb^*}_{i}\cdot \Delta t
\end{equation}

\hfill

\noindent(a-4) Compute nodal intermediate water velocity\\

The nodal intermediate water velocity, $\mathbf{v}^{w^*}_i$, is computed with the water momentum balance equation at the nodes in the two-phase material body:
\begin{equation}
    \mathbf{v}^{w,tb^*}=\frac{\mathbf{b}^{w,tb}_i+ \mathbf{t}^{w,tb}_i+ \mathbf{f}^{w,tb}_i}{m^{w,tb}_i\cdot \frac{ng}{k}} - m^{w,tb}_i \mathbf{a}^{s,tb*}_i
    \label{eq:int_water_acc}
\end{equation}
Where $\mathbf{b}^w_i$ is the nodal water body force, $\mathbf{t}^w_i$ is the nodal water traction, $\mathbf{f}^w_i$ is the nodal wter internal force, and $m^w_i$ is the nodal water mass.

\subsubsection{Pore pressure update (step B)}
The nodal pore pressure increment is implicitly computed globally with the mixture mass balance equation (equation \ref{eq:acc_chorin_mix_mass}):
\begin{equation}
    L(p^{t+1}-p^t)= -A \mathbf{v}^{s^*} - B\mathbf{v}^{wR^*}
\end{equation}
Where, $L$ is the global stiffness matrix, and $A$ and $B$ are the divergence matrices. The global stiffness matrix and divergence matrices are constructed by mapping the element stiffness matrix and divergence matrices. The component of element stiffness matrix and divergence matrices in an element can be computed as:

\begin{equation}
    L_{ij}=\sum_{k=1}^{n_p}[\frac{k}{\rho^wg}(\frac{\rho^w}{\rho}-n)-\frac{\Delta t}{\rho}]\nabla N_i(\mathbf{X}_k)\nabla N_j(\mathbf{X}_k)
\end{equation}
\begin{equation}
    A_{ij}=\sum_{k=1}^{n_p}\nabla N_i(\mathbf{X}_k)\nabla N_j(\mathbf{X}_k)
\end{equation}
\begin{equation}
    B_{ij}=\sum_{k=1}^{n_p}n\nabla N_i(\mathbf{X}_k)\nabla N_j(\mathbf{X}_k)
\end{equation}
where $n_p$ is the number of particles in an element, $N(x)$ is the spatial interpolation function, and subscripts $i$ and $j$ are nodes that belong to the element.

\subsubsection{Final velocity and acceleration (step C)}
\noindent(c-1) Compute final velocities of the two-phase material body\\

The intermediate solid acceleration and velocity are corrected to the final acceleration and velocity with the updated pore pressure increment. The nodal final solid acceleration of the two-phase material body is:
\begin{equation}
    \mathbf{a}^{s.tb}_{i,final}=\frac{\mathbf{b}^{mix,tb}_{i}+\mathbf{t}^{mix,tb}_i+\mathbf{f}^{mix,tb}_i-\nabla(p^{t+1}-p^t)_i}{m^{mix,tb}_i}
    \label{eq:final_tp_solid_acc}
\end{equation}
Where the subscript $final$ denotes the final kinematic parameters at the end of the time step, and the last term of the numerator, $\nabla(p^{t+1}-p^t)_i$, is the nodal gradient of pressure increment. The nodal final solid velocity is updated:
\begin{equation}
    \mathbf{v}^{s,tb}_{i,final}=\mathbf{v}^{s,tb}_{i,initial}+\Delta t \cdot \mathbf{a}^{s,tb}_{i,final}
    \label{eq:final_tp_solid_vel}
\end{equation}

\hfill

\noindent(c-2) Compute final rigid body acceleration

The intermediate rigid body acceleration is adjusted to the final rigid body acceleration with the updated pore pressure increment:

\begin{equation}
     \mathbf{a}^{rb}_{final}=\frac{\splitfrac{\sum_{i=1}^{n_{rb}} (\mathbf{b}^{mix,rb}_i+\mathbf{t}^{mix,rb}_i)+\sum_{i=1}^{n_{contact}} (\mathbf{b}^{mix,tb}_i+\mathbf{t}^{mix,tb}_i)}{+\sum_{i=1}^{n_{contact}} \mathbf{f}^{mix,tb}_i-\sum_{i=1}^{n_{contact}}\nabla(p^{t+1}-t^t)}}{\sum_{i=1}^{n_{contact}} m^{mix,CoM}_i} 
\end{equation}
The final rigid body acceleration is set to the nodal final solid acceleration at the nodes on the rigid body:
\begin{equation}
    \mathbf{a}^{s,rb}_{i,final}=\mathbf{a}^{rb}_{i}
\end{equation}

\hfill

\noindent(c-3) Apply contact algorithm at the interface

To apply the contact algorithm at the interface, the nodal final center-of-mass acceleration is set to the nodal final rigid body acceleration at the nodes at the interface. 
\begin{equation}
    \mathbf{a}^{s,CoM}_{i,final}=\mathbf{a}^{s,rb}_{i,final}
\end{equation}
The nodal final center-of-mass solid velocity at the interface is updated:
\begin{equation}
    \mathbf{v}^{s,CoM}_{i,final}=\mathbf{v}^{s,CoM}_{i,initial}+\Delta t \cdot \mathbf{a}^{s,CoM}_{i,final}
\end{equation}
The contact algorithm is applied at the interface nodes. The normal component of nodal final solid velocities of the two-phase body is adjusted to the nodal final center-of-mass solid velocities.
\begin{equation}
    \Delta\mathbf{v}^{s,tb}_{i,normal}=(\mathbf{v}^{s,CoM}_{i,normal}-\mathbf{v}^{s,tb}_{i,normal})\cdot \mathbf{n}^{tb}_i
\end{equation}
\begin{equation} 
\mathbf{v}^{s,tb}_{i,final}=\mathbf{v}^{s,tb}_{i,final}+\Delta\mathbf{v}^{s,tb}_{i,normal}\cdot \mathbf{n}^{tb}_i
\end{equation}
Lastly, the nodal final kinematics are mapped back to particles, and particle velocity, displacement, location, volume, and porosity are updated.
\subsection{Interface detection process}

\begin{figure}[h]
    \centering
    \includegraphics[width=0.6\textwidth]{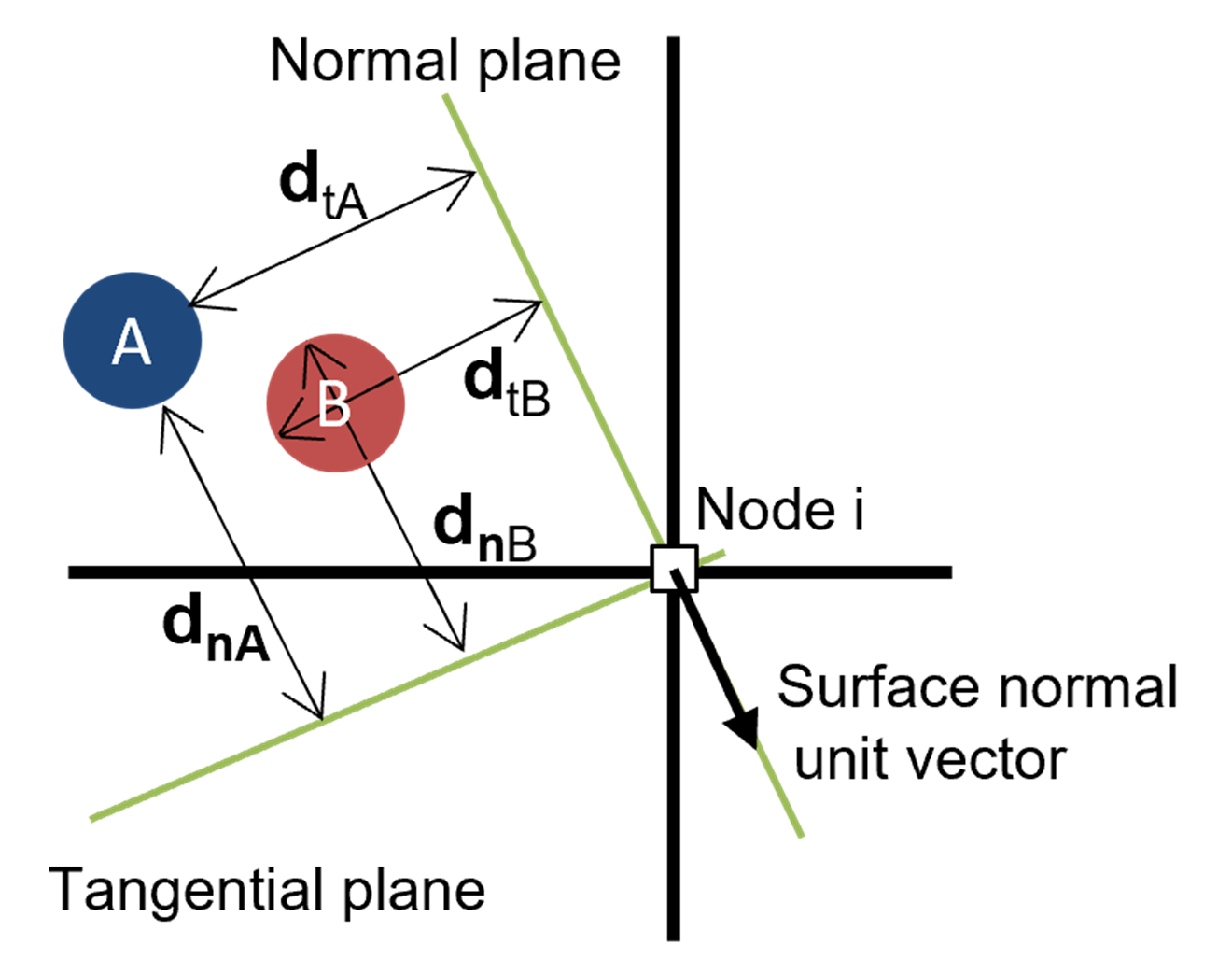}
    \caption{Calculation of edge-to-node distances in normal and tangential directions}
    \label{fig:contact_detection1}
\end{figure}

In order to accurately determine if materials are in contact, we use particle location and sizes to calculate the edge-to-edge distance between particles. The distance between a particle center to a node is computed and projected in the surface normal and tangential directions using a surface normal unit vector at the node. Particle sizes in the normal and tangential directions are updated and utilized to calculate the distance between the particle edge to the node.
\begin{equation}
    d_n=(\mathbf{X}_p-\mathbf{X}_i)\cdot\hat{\mathbf{n}}\pm R_{pn}
\end{equation}
\begin{equation}
    d_t=(\mathbf{X}_p-\mathbf{X}_i)\cdot\hat{\mathbf{t}}\pm R_{pt}
\end{equation}
Where $d_n$ and $d_t$ are the edge-to-node distances in the normal and tangential directions,respectively. $\mathbf{n}$ and $\mathbf{t}$ are the surface normal and tangential unit vectors computed at the node, $\mathbf{X}_p$ is the particle center coordinate, $\mathbf{X}_i$ is the nodal coordinate, and $R_{pn}$ and $R_{pt}$ are deformed particle sizes in the normal and tangential directions. Figure \ref{fig:contact_detection1} illustrates the computation of distances from the particle edge to the node in the normal and tangential directions with a given surface normal unit vector.

\begin{figure}[h]
    \centering
    \includegraphics[width=1\textwidth]{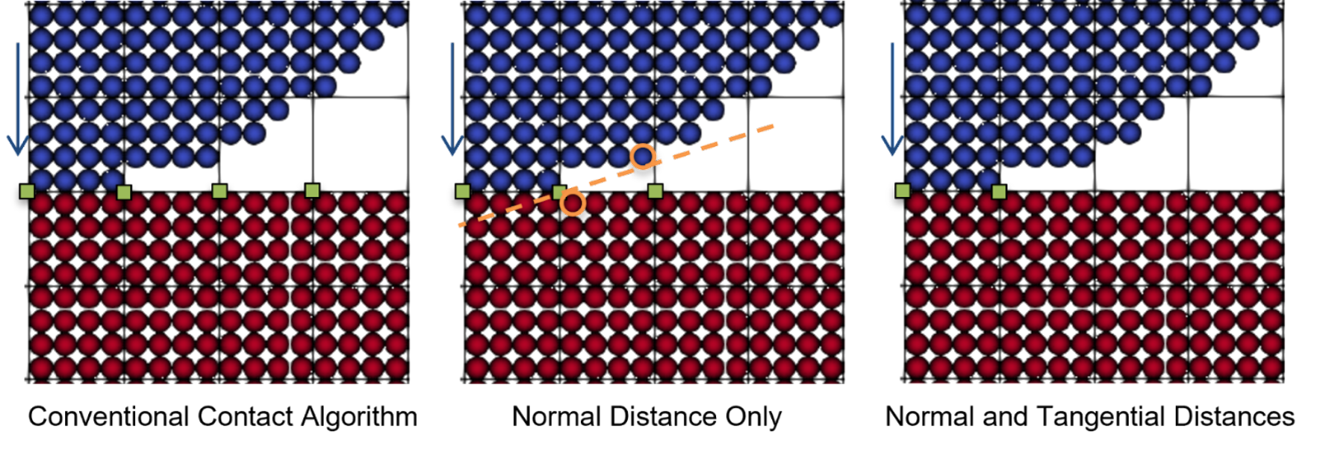}
    \caption{Comparison of contact detection algorithms}
    \label{fig:contact_detection2}
\end{figure}

When particles of different materials approach each other, edge-to-edge distances between particles are updated to detect material contact at nodes. It is important to carefully consider the location of a material body and its particles to determine whether the particle size should be subtracted or added to the center-to-node distances when computing the edge-to-edge distances among particles. For instance, as shown in figure \ref{fig:contact_detection1}, particle size should be added to the center-to-node distance for particle B to correctly compute the edge-to-edge distance between particle A and particle B. For material bodies with simple geometry, using only normal distances can accurately detect the interface. However, for complex geometries with different normal unit vectors at each node, both normal and tangential distances should be used for accurate contact detection process. Figure \ref{fig:contact_detection2} shows the difference in contact node detection results among the conventional contact algorithm, using only normal edge-to-edge distances, and using both the normal and tangential edge-to-edge distances. The contact nodes are marked with a green square. Using both the normal and tangential distances yields the most accurate result. To compute the particle size in the normal and tangential directions, we use the elliptical particle method proposed by \citep{nairn2020new}.

\section{Model Verification}
We evaluate the proposed semi-implicit contact MPM algorithm with three models: (1) one-dimensional consolidation, (2) dynamic impact test, and (3) strip footing. We validate the contact MPM algorithm by comparing the simulation results with those from conventional FEM, explicit MPM, or analytical solutions. For the MPM simulations, we use two-dimensional quadrilateral background elements with linear interpolation functions for displacement computation. The solid bodies are modeled as rigid and the interface between the soil and structure is frictionless. The soil is modeled as porous media fully saturated with water, and the permeability of water is assumed to be constant throughout the simulations.

\subsection{One-dimensional consolidation}
The purpose of the one-dimensional consolidation simulation is to evaluate the accuracy of the semi-implicit contact MPM algorithm, with a focus on minimizing artificial oscillations under instantaneous loading conditions. We analyze pore pressure dissipation and effective stress generation over time along with the soil depth, and compare our simulation results with those of an explicit MPM and Terzaghi's one-dimensional consolidation solution to verify our model.

\begin{figure}[h]
    \centering
    \includegraphics[width=0.4\textwidth]{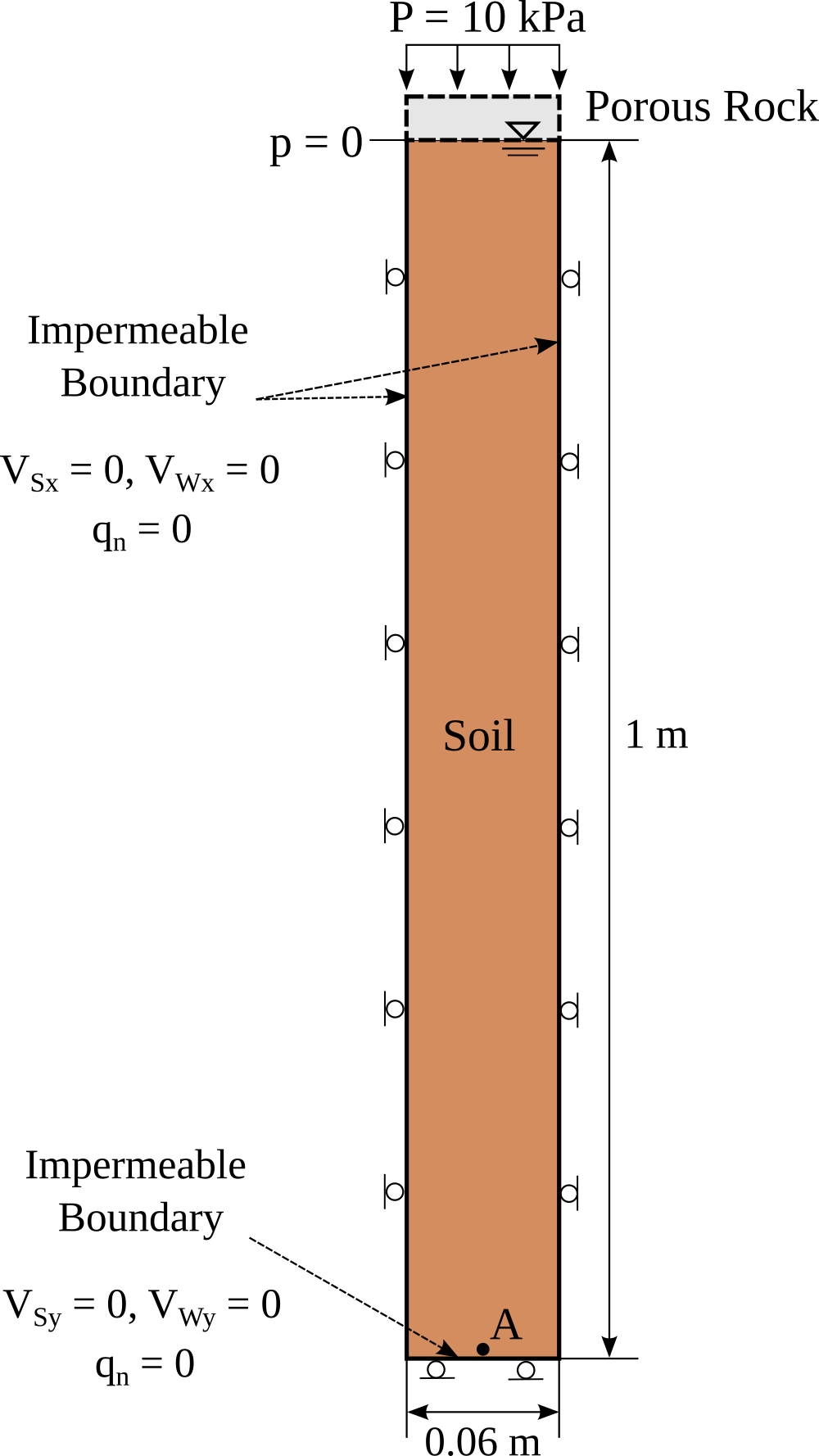}
    \caption{Geometry and boundary conditions of the one-dimensional consolidation simulation}
    \label{fig:consol1}
\end{figure}

The simulation geometry and boundary conditions are depicted in Figure \ref{fig:consol1}. The boundaries at the bottom and sides of the saturated soil column are impervious, and we apply a free-drainage boundary condition at the interface between the soil surface and the loading cap. The soil column is loaded with a instantaneous vertical compressive traction of 10 kPa at the top of the loading cap. The load is transferred to the soil column through contact between the loading cap and the soil surface. Initial pore pressure of 10 kPa and zero initial effective stress are assigned to the soil column. The effect of gravity is not considered.

The soil is isotropic poroelastic, and the loading cap is rigid. The material property of the soil is presented in Table \ref{tab:consol_material}. The total simulation time is 2 seconds with the time step size of 1E-4 seconds. We use quadrilateral elements with the size of 0.02 m by 0.02 m, and 4 particles are located in a cell. For the explicit MPM simulation, we use the same geometry, elements, loading and boundary conditions, and material properties. We use the time step size of 5E-06 seconds to satisfy the stable time increment criteria for explicit analysis.

\begin{table}[ht]
    \centering
    \begin{tabular}{c|c}
    \hline
         Density of solid grain & 2,700 $kg/m^3$ \\
         Porosity & 0.3 \\
         
         Permeability & 1E-3 $m/s$ \\
         
         Young's modulus & 10 $MPa$ \\
         Poisson's ratio & 0.2\\
    \hline
    \end{tabular}
    \caption{Material properties of the soil column}
    \label{tab:consol_material}
\end{table}

\begin{figure}[h]
    \centering
    \includegraphics[width=0.8\textwidth]{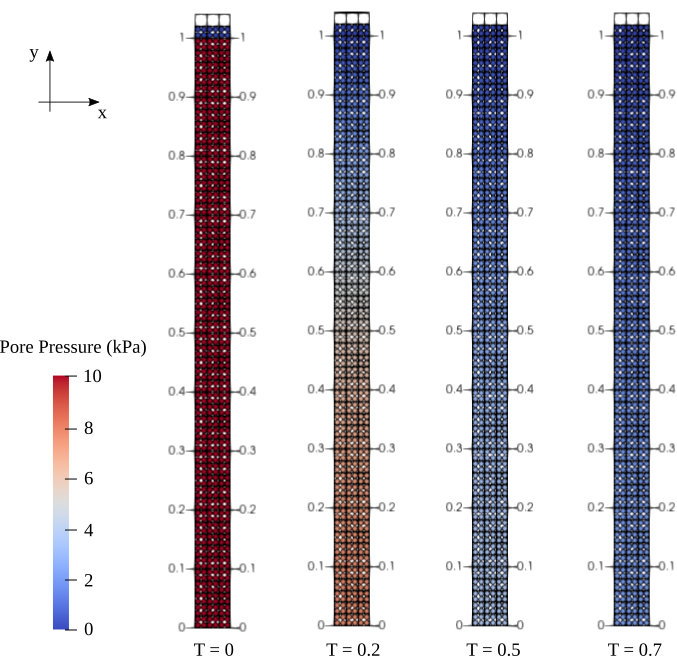}
    \caption{Dissipation of pore pressure with time}
    \label{fig:consol2}
\end{figure}

\begin{figure}[h]
    \centering
    \includegraphics[width=0.8\textwidth]{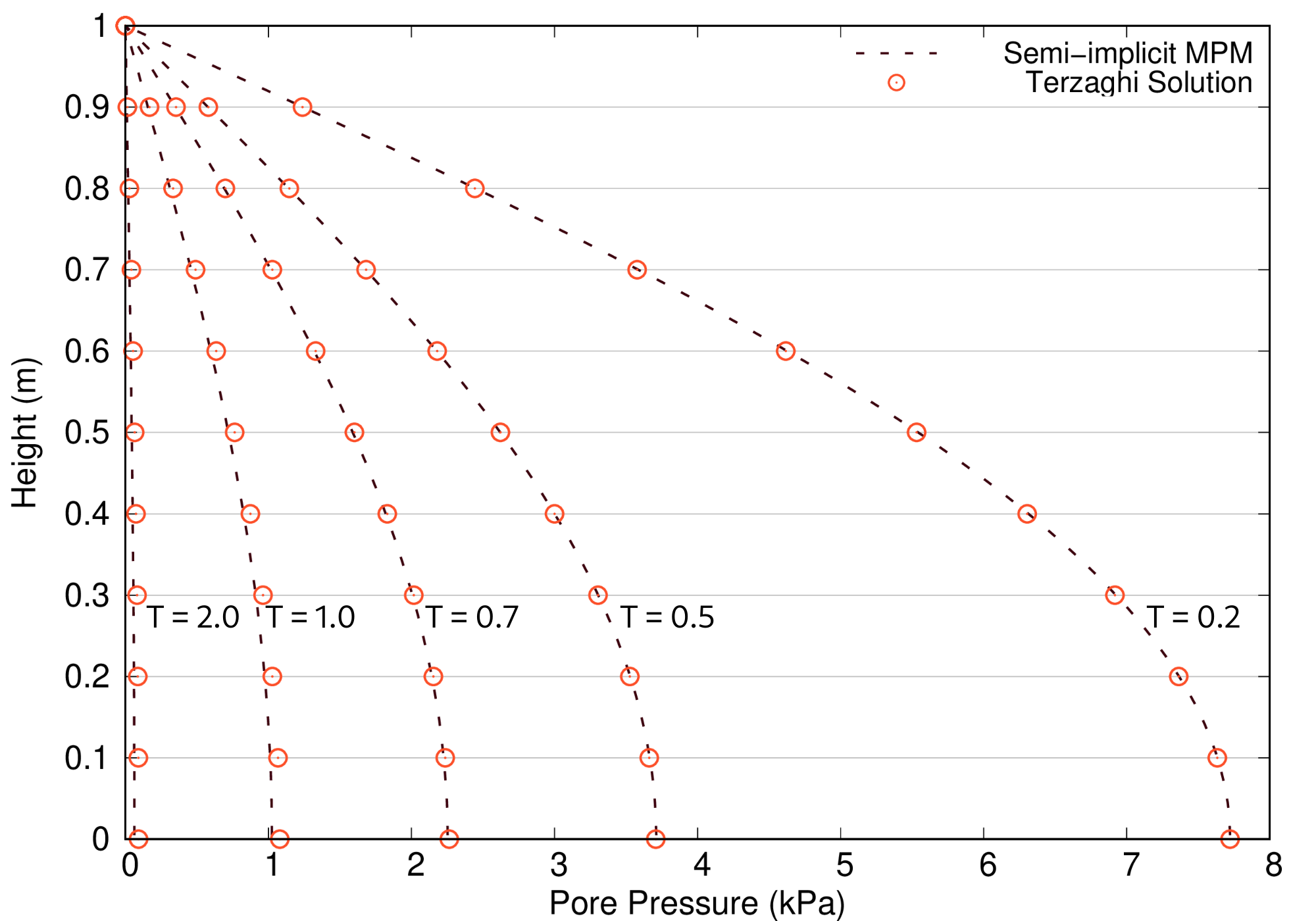}
    \caption{Change in pore pressure distribution with soil column height at different times}
    \label{fig:consol3}
\end{figure}

We present the dissipation of the initial pore pressure from the consolidation simulation in Figure \ref{fig:consol2} and Figure \ref{fig:consol3}. Figure \ref{fig:consol2} shows the pore pressure distribution with time from the semi-implicit contact MPM simulation. We plot the pore pressure with the height of the soil column at different consolidation times in Figure \ref{fig:consol3}. To evaluate the accuracy of the algorithm in the steady state, the results are compared with Terzaghi’s analytical solution, and we confirm that the pore pressure from the semi-implicit MPM simulation agrees well with Terzaghi’s solution.

Figure 9 shows the changes in pore pressure and vertical effective stress at the bottom of the soil column (point A) with time. To consider the dynamic effect from the instantaneous loading, the results from the semi-implicit algorithm are compared to the results of an explicit simulation with the u-p formulation. The plot indicates that the semi-implicit time integration scheme is far more stable than the explicit time integration scheme even with the much larger time step size. The semi-implicit contact MPM significantly reduces the artificial oscillation in pore pressure at the early stage of the simulation. In addition, both pore pressure and effective stress correspond to Terzaghi’s solution at the steady state. The results show the advanced features of the semi-implicit MPM algorithm over the explicit MPM. We also verify the capability of the algorithm to accurately simulate the generation and dissipation of pore pressure and effective stress.

\begin{figure}[htbp]
    \centering
    \begin{subfigure}[b]{0.8\textwidth}
    \includegraphics[width=\textwidth]{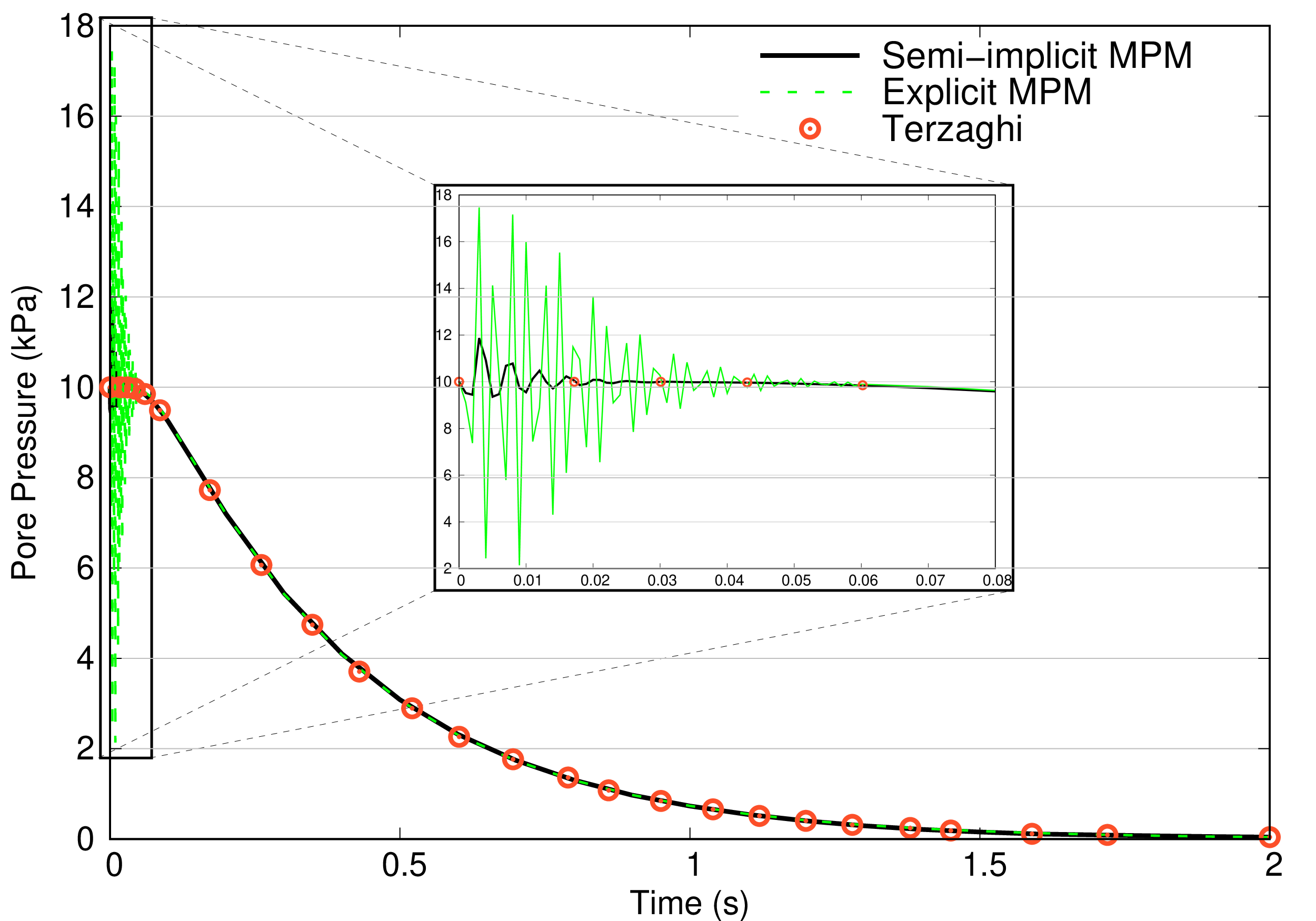}
    \caption{Pore pressure}
    \end{subfigure}\\
    \hfill
    
    \begin{subfigure}[b]{0.8\textwidth}
    \includegraphics[width=\textwidth]{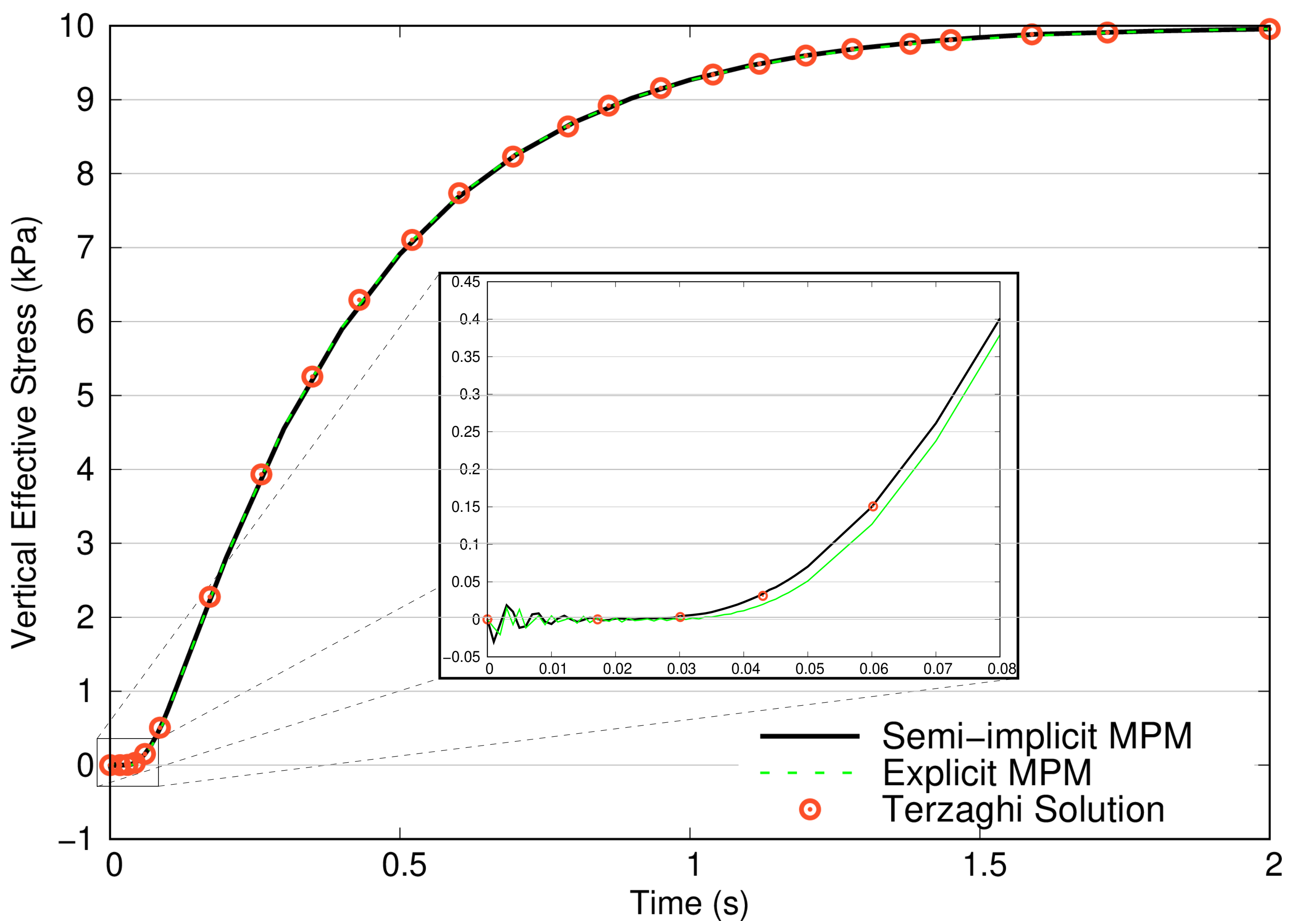}
    \caption{Vertical effective stress}
    \end{subfigure}
    \caption{Change in pore pressure and vertical effective stress at bottom of the soil column (point A) with time}
    \label{fig:ex1_2}
\end{figure}

\subsection{Dynamic impact test}
\begin{figure}[h]
    \centering
    \includegraphics[width=1\textwidth]{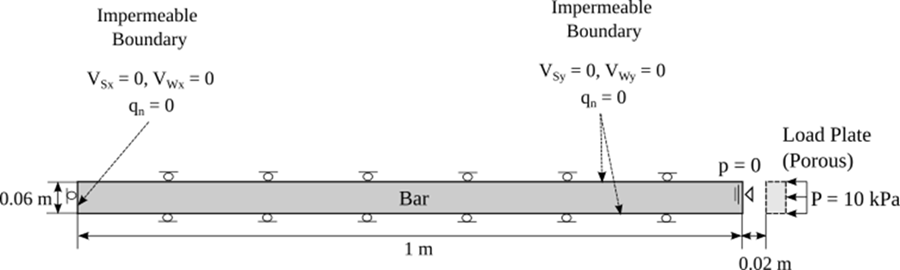}
    \caption{Geometry and boundary conditions of the on-dimensional impact test simulation}
    \label{fig:impact1}
\end{figure}

A simulation of an impact between a moving block and a one-dimensional bar was conducted. The block is a rigid material, and the bar is modeled as a two-phase elastic material, fully saturated with water. The purpose of the one-dimensional dynamic impact simulation is to evaluate the accuracy and stability of the semi-implicit contact MPM algorithm in modeling dynamic problems. The accuracy of the contact detection algorithm is also evaluated. Figure \ref{fig:impact1} shows the geometry and boundary conditions of the model. We model the bar as rectangular with a height of 0.06 m and a width of 1 m.. The water and solid velocities are restricted in the normal direction of each boundary at the bottom, top, and left edge of the bar. We apply the free-drainage condition at the interface between the rigid block and the bar.

\begin{figure}[h]
    \centering
    \includegraphics[width=1\textwidth]{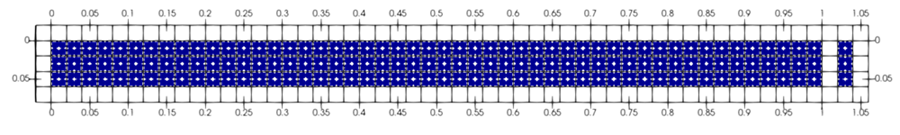}
    \caption{Mesh and particle configuration of the dynamic impact test simulation}
    \label{fig:impact2}
\end{figure}

The instantaneous horizontal load of 10 kPa is applied to the rigid block and maintained throughout the simulation. As the mass of the rigid block is 3.24 kg, the velocity of the rigid block at the time of impact is 2.72 m/s, and the time of impact is 0.0147 seconds. The time step size is 2E-5 seconds, and the simulation lasted until it reaches the steady state. We use quadrilateral elements with the size of 0.02 m by 0.02 m, and 4 particles are located in a cell. The mesh and particle configuration are shown in Figure \ref{fig:impact2}, and the material properties of the bar and the rigid block are summarized in table \ref{tab:impact_material}.

\begin{table}[ht]
    \centering
    \begin{tabular}{c|c|c}
    \hline
    &Density of solid grain& 2,700$kg/m^3$ \\
    &porosity & 0.3\\
    Soil&Permeability & 1E-3 $m/s$\\
    &Young's modulus&2 $MPa$\\
    &Poisson's ratio&0.2\\
    \hline
    Rigid block&Mass&3.24 $kg$\\     
    \hline
    \end{tabular}
    \caption{Material properties of the soil and rigid block}
    \label{tab:impact_material}
\end{table}

\begin{figure}[h]
    \centering
    \includegraphics[width=1\textwidth]{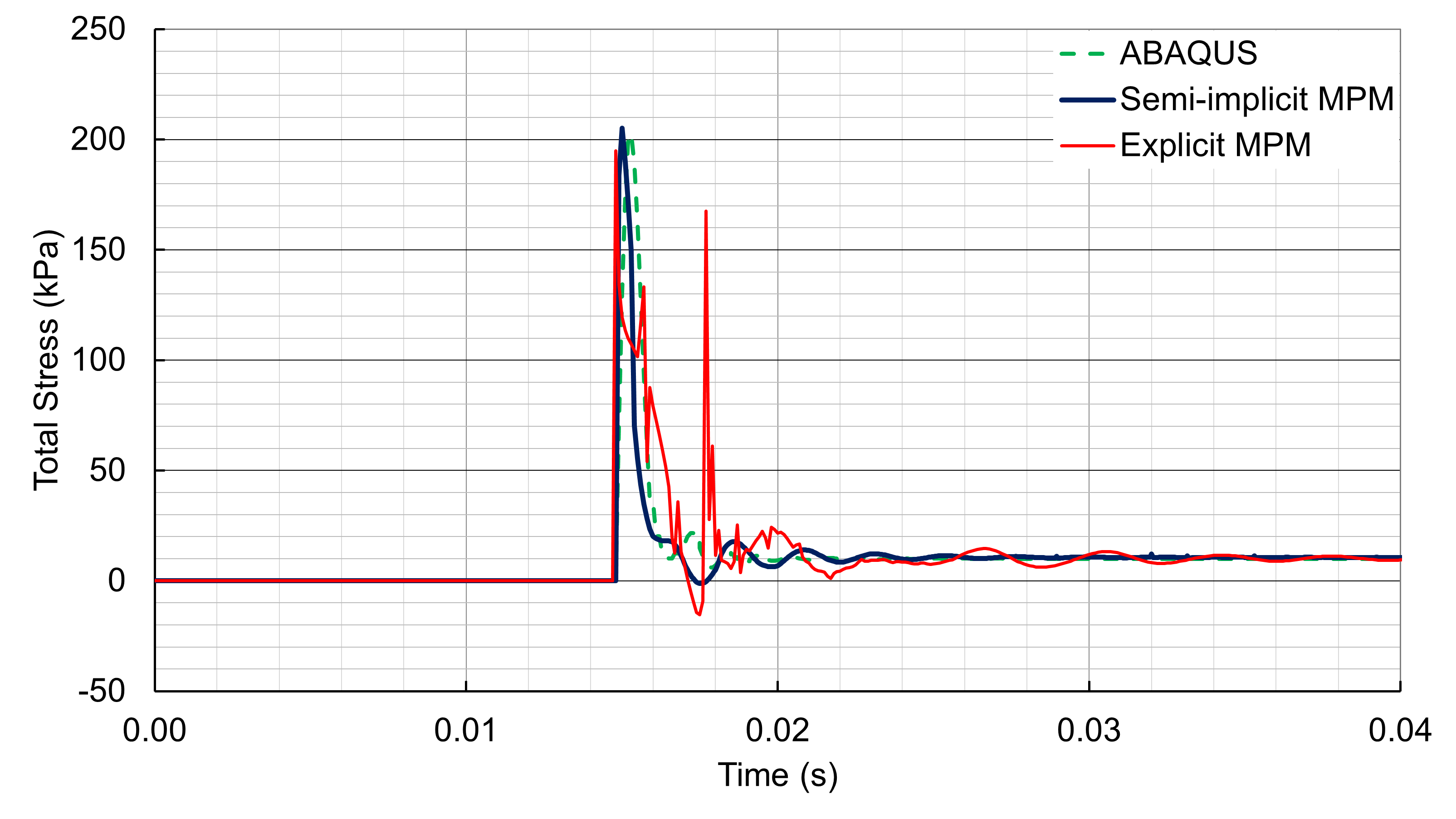}
    \caption{Comparison of total stresses at the interface}
    \label{fig:impact3}
\end{figure}

We monitor the total stress and pore pressure waves generated in the bar from the impact, and the results are compared with an explicit MPM simulation and an ABAQUS simulation. The identical geometry, boundary, and loading conditions are applied to the explicit and ABAQUS simulations. The step size of 1E-6 seconds is used for the explicit MPM simulation. Figure \ref{fig:impact3} shows the total stress wave at the time of the impact measured at the interface between the bar and rigid block. All simulations yield similar peak stress values from the first wave, but the explicit MPM simulation has a severe numerical oscillation problem and yields a significant error for the first several waves. The explicit simulation also takes longer to reach the steady state. It is also noted that the proposed contact detection algorithm yields an accurate time of impact as 0.0147 seconds.

\begin{figure}[h]
    \centering
    \includegraphics[width=1\textwidth]{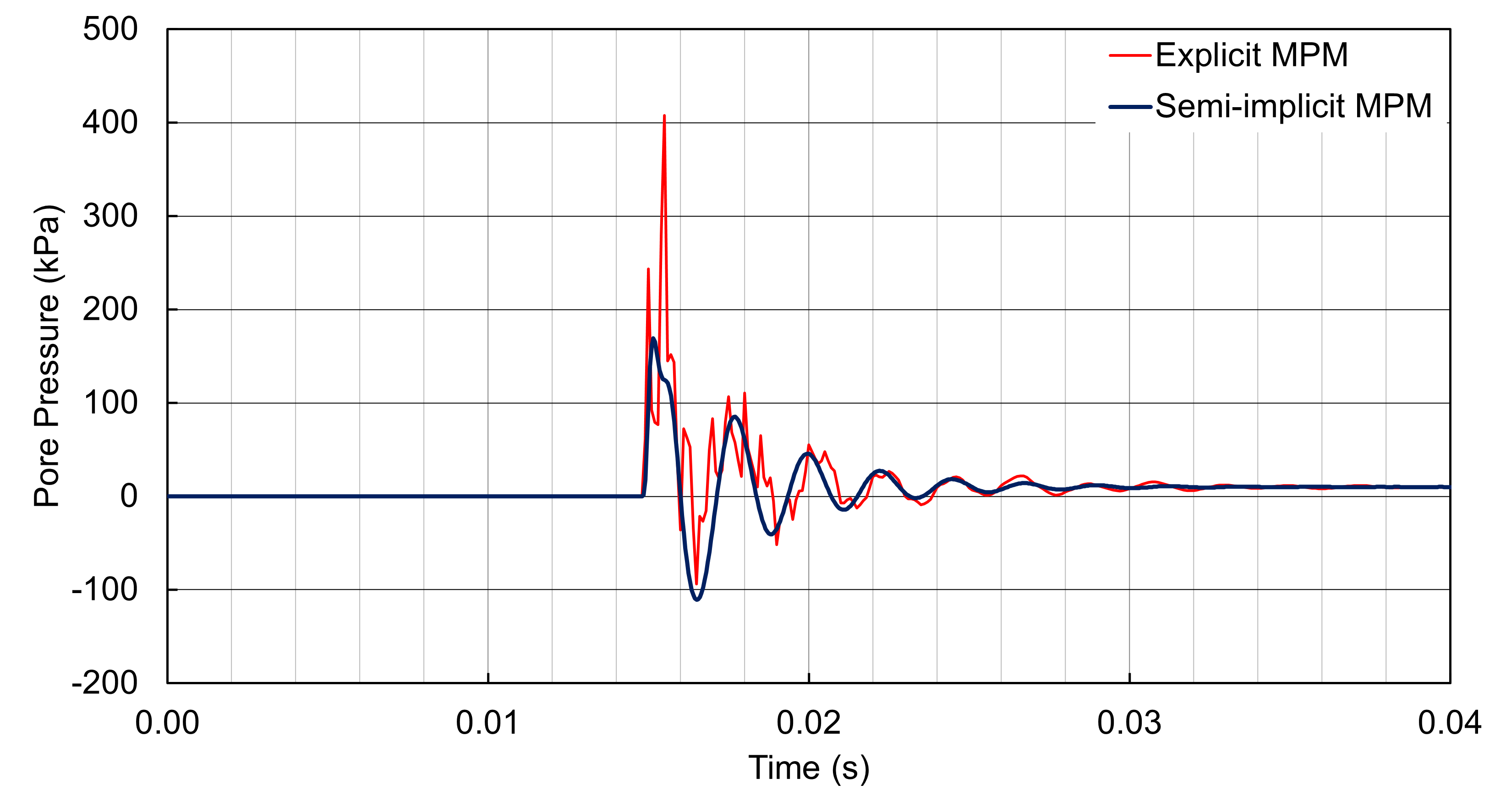}
    \caption{Pore pressure wave at the center of the two-phase bar}
    \label{fig:impact4}
\end{figure}

 Figure \ref{fig:impact4} shows the pore pressure wave measured at the center of the bar from the semi-implicit and explicit simulations. The peaks and wave velocities from the semi-implicit and explicit simulations generally agree well with each other, while the explicit simulation yields a significant error in the peak pressure. Notably, the semi-implicit contact algorithm does not suffer from pressure oscillations and severe errors in solving two-phase dynamic problems with a larger time step size.

\subsection{Strip foundation simulation}
\begin{figure}[h]
    \centering
    \includegraphics[width=1\textwidth]{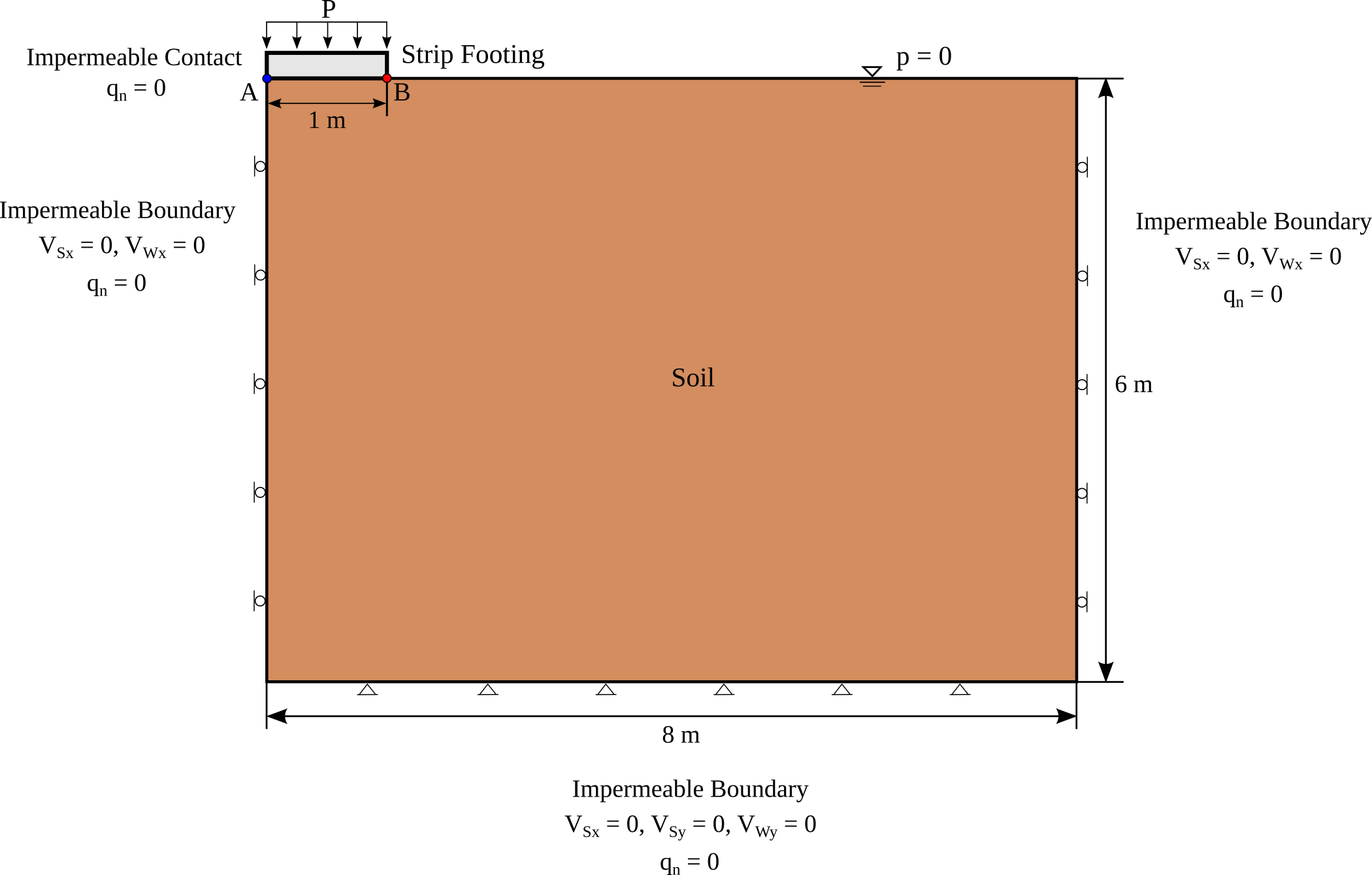}
    \caption{Geometry and boundary conditions of strip foundation simulation}
    \label{fig:FD1}
\end{figure}

The purpose of the strip foundation simulation is to examine the capability of the semi-implicit contact algorithm in modeling: (1) the response of the saturated soil near an undrained soil-structure interface, and (2) two-dimensional pore pressure generation and dissipation, and to compare the MPM simulation to the conventional FEM simulation. We simulate the strip foundation with two-dimensional plain-strain elements. We model half of the foundation and the soil, taking advantage of the symmetric geometry. Figure \ref{fig:FD1} shows the geometry and boundary conditions. The width of the foundation is 1 m, and the width and depth of the soil are 8 m and 6 m. The bottom and the sides of the soil are impermeable, and we apply the free-drainage condition at the soil surface including point B. The interface between the soil and the strip foundation is impermeable. At the bottom of the background grid, we apply the fully fixed boundary conditions, while only horizontal velocities are restricted at the sides.

\begin{figure}[h]
    \centering
    \includegraphics[width=0.6\textwidth]{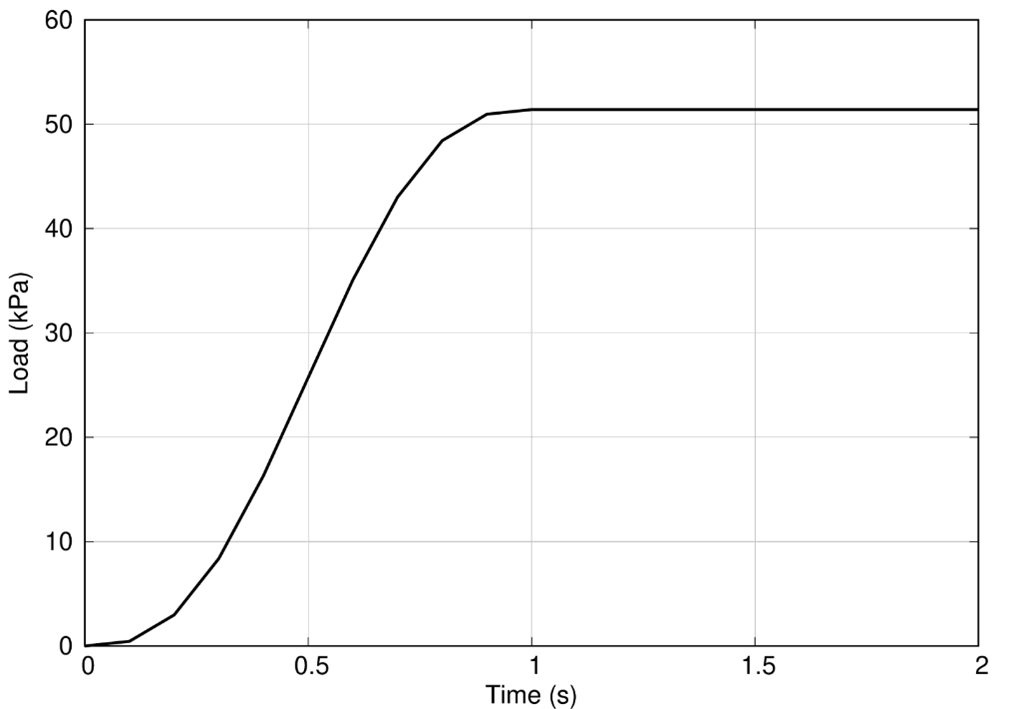}
    \caption{Load history of strip foundation simulation}
    \label{fig:FD2}
\end{figure}

Tresca material model was used to simulate the soil behavior. We apply a 51.4 kPa of vertical compressive traction, which is Prandtl’s bearing capacity solution (Prandtl, 1920), at the top of the foundation. This traction load increases from 0 kPa to 51.4 kPa gradually in one second and remains constant after. The load history is shown in Figure \ref{fig:FD2}. The total simulation time is two seconds, and the time step size is 1E-4 seconds. We use quadrilateral elements with a size of 0.2 m by 0.2 m, and nine particles are located in a cell.
We monitor the generation of pore pressure and effective stresses near the undrained interface, the displacement of the foundation, and the effective stress concentration and failure of soil at the edge of the interface (point B). The material properties used in the MPM simulation are presented in Table \ref{tab:FD_material}. 

\begin{table}[ht]
    \centering
    \begin{tabular}{c|c}
    \hline
         Density of solid grain & 2,150 $kg/m^3$ \\
         Porosity & 0.2 \\
         
         Permeability & 5E-3 $m/s$ \\
         
         Young's modulus & 1 $MPa$ \\
         Poisson's ratio & 0.2\\
         Undrained Shear Strength & 10 $kPa$ \\
    \hline
    \end{tabular}
    \caption{Material properties of the soil}
    \label{tab:FD_material}
\end{table}

We conduct the finite element simulation with ABAQUS 2018. We use identical geometry, boundary and loading conditions, and material properties with the MPM simulation. The strip foundation is modeled as a rigid body, and we apply the penalty contact at the interface between the soil and the foundation. Two-dimensional plane strain square elements with the size of 0.1 m by 0.1 m are used, and the maximum time step size is set to 1E-3 seconds. The ABAQUS simulation failed at 0.72 seconds due to a convergence problem with the attempted time step size of 1E-9 seconds.

\begin{figure}[h]
    \centering
    \includegraphics[width=1\textwidth]{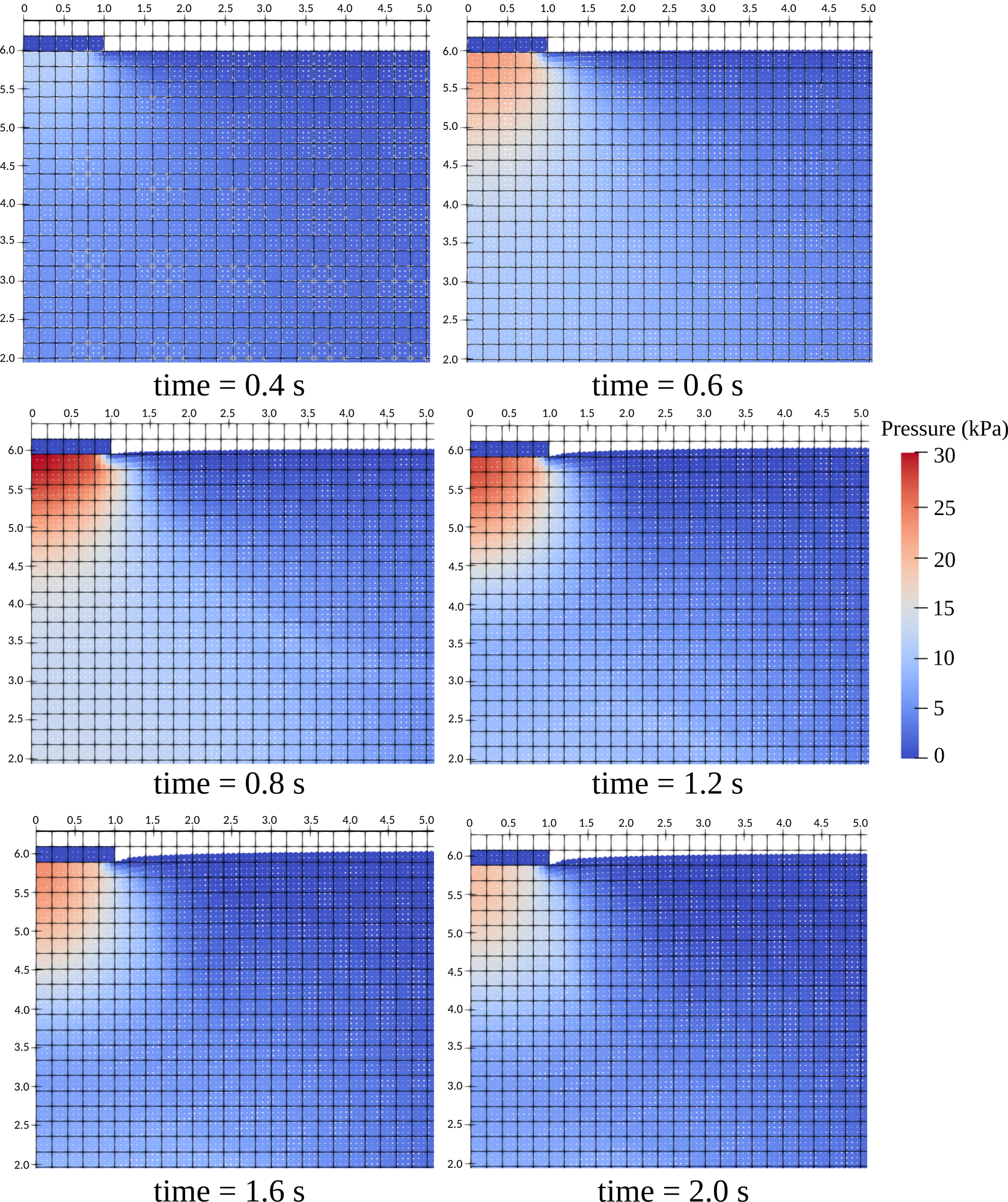}
    \caption{Comparison of pore pressure changes with time at point A from MPM and ABAQUS simulations}
    \label{fig:FD3}
\end{figure}

\begin{figure}[h]
    \centering
    \includegraphics[width=0.8\textwidth]{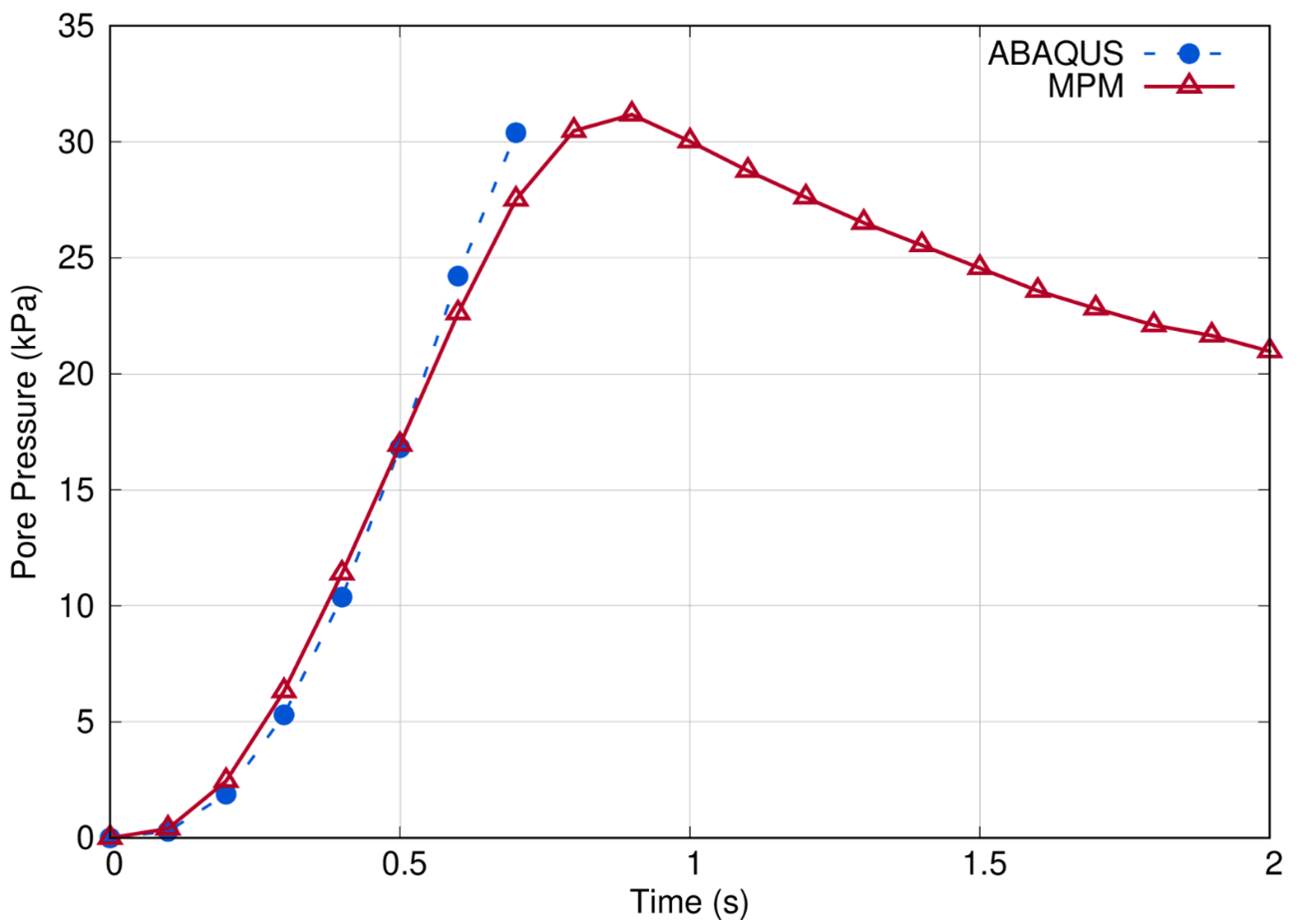}
    \caption{Comparison of the pore pressure changes with time at apoint A from MPM and ABAQUS simulations}
    \label{fig:FD4}
\end{figure}

Figure \ref{fig:FD3} presents the pore pressure distributions from the MPM simulation with time, and Figure \ref{fig:FD4} shows the pore pressure at the soil below the center of the foundation (point A) from the MPM and ABAQUS simulations. The generation and dissipation of the pore pressure with the undrained interface are captured accurately with the MPM contact algorithm, and the pore pressure from the MPM simulation corresponds with the ABAQUS simulation.

 \begin{figure}[h]
    \centering
    \includegraphics[width=1\textwidth]{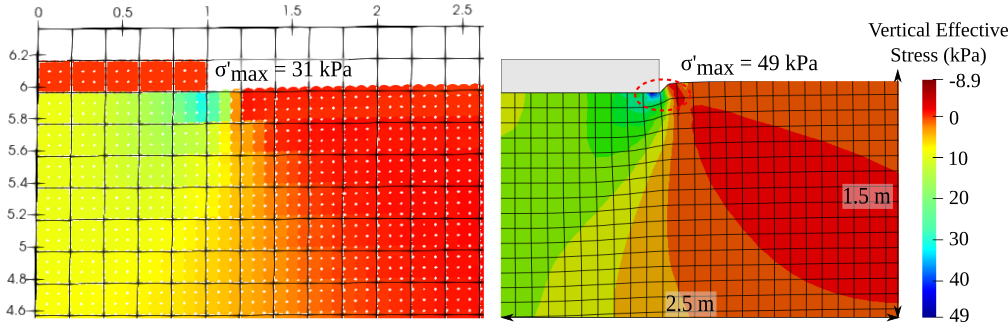}
    \caption{Comparison of vertical effective stress distribution from MPM (left) and ABAQUS (right) simulations at the simulation time of 0.7 seconds}
    \label{fig:FD5}
\end{figure}

\begin{figure}[h]
    \centering
    \includegraphics[width=0.8\textwidth]{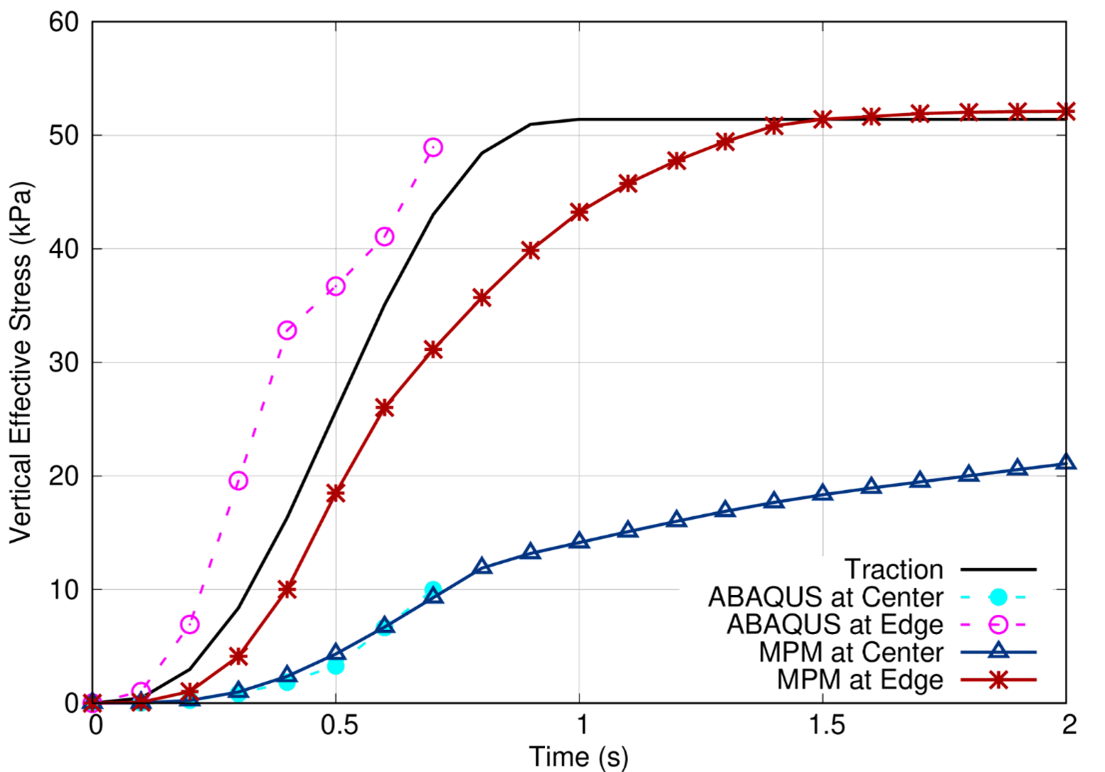}
    \caption{Change of vertical effective stress at point A and point B from MPM and ABAQUS simulations}
    \label{fig:FD6}
\end{figure}

Figure \ref{fig:FD5} shows the vertical effective stress distribution from the MPM and ABAQUS simulations when the simulation time is 0.7 seconds, at which the vertical traction applied at the footing is 43.02 kPa. The ABAQUS simulation exhibits an excessive stress concentration at the soil surface beneath the foundation's edge (point B), leading to significant distortion of the soil element as it deforms. Figure \ref{fig:FD6} presents the vertical effective stresses measured at the soil below the center (point A) and edge (point B) of the foundation from the MPM and ABAQUS simulations. Both the MPM and ABAQUS simulations yield identical results at point A, but the ABAQUS simulation significantly overestimates the vertical effective stress at point B, as the vertical effective stress exceeds the traction load.

\begin{figure}[h]
    \centering
    \includegraphics[width=0.85\textwidth]{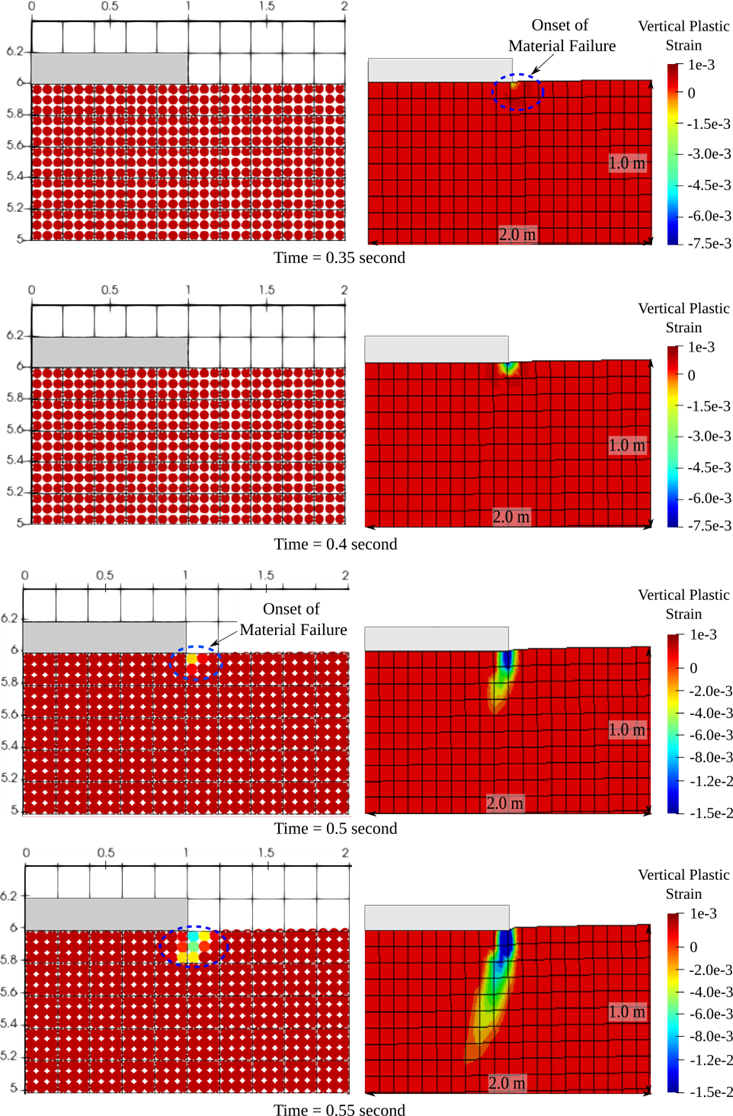}
    \caption{Development of vertical plastic strain over time from MPM (left) and ABAQUS (right) simulations}
    \label{fig:FD7}
\end{figure}

\begin{figure}[h]
    \centering
    \includegraphics[width=0.8\textwidth]{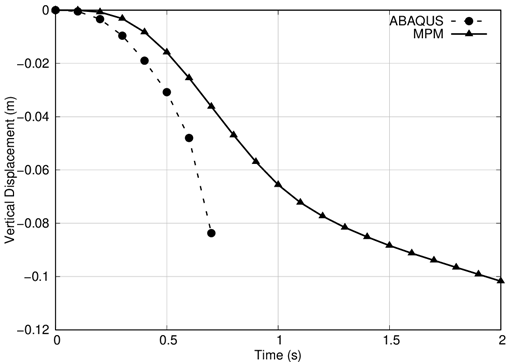}
    \caption{Comparison of vertical displacement of foundation with time from MPM and ABAQUS simulations}
    \label{fig:FD8}
\end{figure}

The excessive effective stress from the singularity problem results in the failure of the soil material in the earlier stage. Therefore, the ABAQUS simulation overestimates the overall foundation displacement, thus, the response of the saturated soil is not accurately simulated. Figure \ref{fig:FD7} illustrates the vertical plastic strain distribution from the MPM and ABAQUS simulations with time. The plastic strain distribution shows the part at which the material failure has been initiated. The MPM simulation shows the onset of plasticity between 0.5 and 0.55 seconds, with traction between 25.7 kPa and 30.49 kPa, close to the lower-bound capacity for initial plasticity of 31.4 kPa. In contrast, the ABAQUS simulation exhibits the onset of plasticity at 0.35 seconds, with a traction of 12.09 kPa, and also yields a larger area of failed soil at the same simulation time. Figure \ref{fig:FD8} shows that the foundation displacement is more than twice as large in the ABAQUS simulation due to the early material from the stress singularity.

The two-dimensional strip foundation simulation demonstrates the efficacy of the semi-implicit contact algorithm in handling soil-structure interaction problems with an undrained interface. The semi-implicit contact MPM can accurately simulate the generation and dissipation of pore pressure and effective stresses while avoiding the conventional finite element method's interface problems, such as stress singularity and element distortions.

\section{Discussions and conclusions}
We develop a semi-implicit two-phase contact algorithm for soil-structure interaction problems in the framework of the MPM. With the introduction of the splitting algorithm and the semi-implicit time increment scheme, the two-phase algorithm inherits significant advantages over explicit dynamic formulations. The algorithm overcomes the artificial pressure oscillation and allows a time increment step independent of fluid compressibility. The rigid body algorithm is introduced to avoid interaction problems between material bodies with significantly different stiffnesses. Lastly, we use the edge-to-edge distances between particles in both the surface normal and tangential directions to accurately identify the nodes at the interface. 

We evaluate the performance of the contact algorithm with practical example simulations. With the one-dimensional consolidation simulation, we verify the accuracy and stability of the algorithm over the explicit time integration scheme. The semi-implicit contact MPM algorithm significantly reduces the artificial oscillation from the explicit MPM even with the larger time step size. The proposed contact algorithm yields identical pore pressure and effective stress distributions with the analytical solution at the steady state and accurately models the consolidation process. The semi-implicit contact algorithm is capable of modeling dynamic problems, and the proposed contact detection algorithm perfectly identifies the contact nodes. The dynamic impact simulation, again, verified the stability of the algorithm with accurate and smooth stress and pressure waves. With the two-dimensional strip foundation simulation, we confirm the capability of the contact algorithm in simulating the pore pressure generation and the two-dimensional consolidation near the undrained interface. The semi-implicit contact MPM simulation has significant advantages over the conventional finite element simulation as the MPM simulation does not suffer from stress singularity and excessive mesh distortion problems.

While the stability of the two-phase contact algorithm has been verified, there is still room for improvement. The current implementation of the splitting method and implicit update of pore pressure requires additional parameters and calculations, resulting in longer simulation times than explicit formulations. Furthermore, using finer elements for more accurate results exponentially increases computation time. To address these issues, implementing advanced parallel computing is highly required. The current algorithm only allows for fully impervious interface or free-drainage conditions at the soil-structure interface. To provide detailed analyses for complex interface problems, it is suggested to implement advanced drainage conditions at the interface. Lastly, to handle highly dynamic loading problems such as liquefaction, it is suggested to develop an algorithm with $u-p-w$ formulations.




\bibliographystyle{elsarticle-num-names}
\bibliography{paper.bib}



\end{document}